\nonstopmode
\nonstopmode
\documentclass[10pt]{amsart}
\usepackage{amssymb}
\usepackage{amscd}
\usepackage{graphicx}
\usepackage{psfrag}
\usepackage{verbatim}
\usepackage[all]{xy}
\newdir{ >}{{}*!/-10pt/@{>}}
\newdir^{ (}{{}*!/-5pt/@^{(}}
\newdir_{ (}{{}*!/-5pt/@_{(}}
\def\undersetbrace#1\to#2{\underbrace{#2}_{#1}}
\def\oversetbrace#1\to#2{\overbrace{#2}^{#1}}
\def\AMSunderset#1\to#2{\underset{#1}{#2}}
\def\AMSoverset#1\to#2{\overset{#1}{#2}}

\newenvironment{proclaim}[1]{\par\medskip\noindent{\bf #1.}\it}{\par\smallskip}
\newenvironment{demo}[1]{\par\smallskip{\bf #1.}}{\par\smallskip}
\newcounter{locnmb}[subsection]
\newcommand{\nmbzero}[1]{\protect\refstepcounter{section}\thesection\label{nmb:#1}}
\newcommand{\nmbdot}[1]{\protect\refstepcounter{subsection}%
\thesubsection\label{nmb:#1}\gdef\lastnmb{#1}}
\newcommand{\nmbcolom}[1]{\protect\refstepcounter{subsubsection}%
\arabic{subsubsection}\label{nmb:\lastnmb.#1}%
\refstepcounter{locnmb}\label{locnmb:\lastnmb.#1}}
\newcommand{\nmbref}[1]{\protect\ref{nmb:#1}}
\newcommand{\nmbvert}[1]{\protect\ref{locnmb:\lastnmb.#1}}
\newcommand{\nmb}[2]{%
\ifx0#1\nmbzero{#2}\else%
\ifx.#1\nmbdot{#2}\else%
\ifx:#1\nmbcolom{#2}\else%
\ifx!#1{\rm(\nmbref{#2})}\else%
\ifx|#1{\rm(\nmbvert{#2})}\else%
\fi\fi\fi\fi\fi%
}
\def\cit#1#2{\ifx#1!\cite{#2}\else#2\fi} 
\def\fig#1#2{\ifx#1!\ref{fig:#2}\else\label{fig:#2}\fi} 
\def\totoc{}             
\def\idx{}               
\def\ign#1{}             
\def\o{\circ}
\def\al{\alpha}
\def\be{\beta}
\def\ga{\gamma}

\def\ep{\varepsilon}

\def\vartheta{\theta}

\def\la{\lambda}

\def\si{\sigma}

\def\ph{\varphi}

\def\Ga{\Gamma}
\def\De{\Delta}

\def\La{\Lambda}

\def\i{^{-1}}
\def\x{\times}
\def\g{\mathfrak g}
\def\Id{\operatorname{Id}}

\def\Fl{\operatorname{Fl}}
\def\p{\partial}
\def\on{\operatorname}
\def\today{\ifcase\month\or
 January\or February\or March\or April\or May\or June\or
 July\or August\or September\or October\or November\or December\fi
 \space\number\day, \number\year}

\newif\ifFIRST\newdimen\MAXright\MAXright0pt
\def\sdynkin{\bgroup\eightpoint\dynkin}
\def\endsdynkin{\enddynkin\egroup}
\def\dynkin{\bgroup\FIRSTtrue\hskip.5em\setbox1\hbox{$\diagup$}%
\setbox2\hbox{$\diagdown$}%
\setbox0\hbox to2\wd1{\hrulefill}%
\setbox3\hbox{$\bullet$}%
\setbox4\hbox{$\circ$}%
\def\root##1{\ifFIRST\setbox5\hbox{$##1$}\ifdim\wd5>1.3em%
\hskip-.5em\hskip.5\wd5\fi\fi\FIRSTfalse%
\hskip-.25em\raise1.5\wd3\hbox to0pt{\hss\hskip.45em$%
\scriptstyle##1$\hss}\copy3\hskip-.25em\setbox6\hbox{$##1$}%
\MAXright\wd6}%
\def\droot##1{\ifFIRST\setbox5\hbox{$##1$}\ifdim\wd5>1.3em%
\hskip-.5em\hskip.5\wd5\fi\fi\FIRSTfalse%
\hskip-.25em\lower1.8\wd3\hbox to0pt{\hss\hskip.45em$%
\scriptstyle##1$\hss}\copy3\hskip-.25em\setbox6\hbox{$##1$}%
\MAXright\wd6}%
\def\rroot##1{\hskip-.25em\copy3\hbox to0pt{\hskip.3em$\scriptstyle##1$\hss}%
\hskip-.25em\setbox6\hbox{\hskip.6em$##1##1$}%
\MAXright\wd6}%
\def\wrroot##1{\hskip-.36em\copy4\hbox to0pt{\hskip.3em$\scriptstyle##1$\hss}%
\hskip-.48em\setbox6\hbox{\hskip.6em$##1##1$}%
\MAXright\wd6}%
\def\wroot##1{\ifFIRST\setbox5\hbox{$##1$}\ifdim\wd5>1.3em%
\hskip-.5em\hskip.5\wd5\fi\fi\FIRSTfalse%
\hskip-.36em\raise1.5\wd3\hbox to0pt{\hss\hskip.6em$%
\scriptstyle##1$\hss}\copy4\hskip-.38em\setbox6\hbox{$##1$}%
\MAXright\wd6}%
\def\wdroot##1{\ifFIRST\setbox5\hbox{$##1$}\ifdim\wd5>1.3em%
\hskip-.5em\hskip.5\wd5\fi\fi\FIRSTfalse%
\hskip-.36em\lower1.8\wd3\hbox to0pt{\hss\hskip.6em$%
\scriptstyle##1$\hss}\copy4\hskip-.38em\setbox6\hbox{$##1$}%
\MAXright\wd6}%
\def\link##1{\raise.22em\copy0\hskip-.5\wd0\hbox to0pt{\hss$##1$\hss}\hskip.5\wd0}%
\def\llink##1{\raise.32em\copy0\hskip-\wd0%
\raise.12em\copy0\hskip-.5\wd0\hbox to0pt{\hss$##1$\hss}\hskip.5\wd0}%
\def\lllink##1{\raise.22pt\copy0\hskip-\wd0\raise.32em\copy0\hskip-\wd0%
\raise.12em\copy0\hskip-.5\wd0\hbox to0pt{\hss$##1$\hss}\hskip.5\wd0}%
\def\rootupright##1{\hbox to0pt{\raise.45em\copy1\hskip-.25em\raise1.3\ht1%
\hbox{\copy3\hskip.3em$\scriptstyle##1$}\hss}%
\setbox6\hbox{\hskip.6em\copy1\copy1$##1##1$}%
\ifdim\MAXright<\wd6\MAXright\wd6\fi}%
\def\rootdownright##1{\hbox to0pt{\raise-.5em\copy2\hskip-.25em\raise-1.35\ht1%
\hbox{\copy3\hskip.3em$\scriptstyle##1$}\hss}\setbox6%
\hbox{\hskip.6em\copy2\copy2$##1##1$}%
\ifdim\MAXright<\wd6\MAXright\wd6\fi}%
\def\rootdown##1{\hbox to0pt{\hskip-.05em\vrule height.25em depth.65em%
\hskip-.25em\raise-.95em\hbox{\copy3\hskip.3em$\scriptstyle##1$}\hss}%
\setbox6\hbox{$##1$}%
\ifdim\MAXright<\wd6\MAXright\wd6\fi}%
\def\dots{\hskip.5em\cdots\hskip.5em}}%
\def\enddynkin{\ifdim\MAXright>1em\hskip.5\MAXright\else\hskip.5em\fi\egroup}%

\def\AMSonly#1{}
\begin{document}
\title{Reflection groups on Riemannian manifolds}
\author[Alekseevski, Kriegl, Losik, Michor]{Dmitri Alekseevsky,
Andreas Kriegl, Mark Losik, Peter W. Michor}
\thanks{P.W.M. and M.L. were supported
by `Fonds zur F\"orderung der wissenschaftlichen
Forschung, Projekt P~14195~MAT'}

\begin{abstract}
We investigate discrete groups $G$ of isometries of a complete 
connected Riemannian manifold $M$ which are generated by reflections, 
in particular those generated by disecting reflections. 
We show that these are Coxeter groups, 
and that the the orbit space $M/G$ is isometric to a Weyl chamber $C$ 
which is a Riemannian manifold with corners and certain angle conditions 
along intersections of faces. We can also reconstruct the manifold and its action from the 
Riemannian chamber and its equipment of isotropy group data along the faces. 
We also discuss these results from the point of view of Riemannian orbifolds.  
\end{abstract}

\address{D.V. Alekseevsky:
Department of Mathematics,
University of Hull,
Cottingham Road,
Hull, HU6 7RX,
England.}
\email{d.v.alekseevsky@maths.hull.ac.uk}
\address{A.\  Kriegl: Institut f\"ur Mathematik, Universit\"at Wien,
Strudlhofgasse 4, A-1090 Wien, Austria.}
\email{Andreas.Kriegl@univie.ac.at}
\address{M. Losik: Saratov State University, ul. Astrakhanskaya, 83,
410026 Saratov, Russia.}
\email{LosikMV@info.sgu.ru}
\address{P.\  W.\  Michor: Institut f\"ur Mathematik, Universit\"at Wien,
Strudlhofgasse 4, A-1090 Wien, Austria; {\it and}:
Erwin Schr\"odinger Institute of Mathematical Physics, Boltzmanngasse
9, A-1090 Wien, Austria.}
\email{Peter.Michor@esi.ac.at}
\def\LaTeXonly{}

\maketitle

\section*{\totoc\nmb0{1}. Introduction }

The aim of this paper is to study the discrete groups $G$
generated by reflections with respect to hypersurfaces (shortly,
reflection groups ) on a Riemannian manifold $M$.  If $M=E^n$ is
the Euclidean space, then the classification of all reflection
groups was given in a fundamental paper by Coxeter \cit!{11}. This
implies also the classification of reflection groups  on the
sphere $S^n$. There are many results about reflection groups in
hyperbolic space, see Vinberg \cit!{32}, \cit!{34}, \cit!{33},
and \cit!{35}, but the complete classification is missing. In all
these cases the appropriate fundamental domain $C$ of a reflection
group $G$ (called Weyl chamber) is a Coxeter polyhedron, i.e., a
convex polyhedron where any two neighbour  walls (codimension 1
faces $F_i$, $F_j$ with codimension 2 intersection)  have angle
$\pi/n_{i,j}$ for $n
_{i,j}\in \mathbb N$. We call this the {\it Coxeter property}. 
Conversely, any Coxeter polyhedron $C$ in a space of
constant curvature $M=S^n, E^n, H^n$ is the fundamental domain of
the reflection group $G$ which is generated by the reflections
$s_i=s_{F_i}$ with respect to the walls $F_i$ of $C$. The group
$G$ is a Coxeter group, i.e., a group with a set
$S=\{s_1,\dots,s_l\}$ of generators, and relations $s_i^2=1$,
$(s_is_j)^{n_{i,j}}=1$  for $n_{i,j}\in \mathbb N\cup\{\infty\}$. In
our case, $n_{i,j}$ is defined by the angle between the walls
$F_i$ and $F_j$ as above.

The manifold $M$ with the action of $G$ can be reconstructed from
the Weyl  chamber  $C$ (which is homeomorphic to the orbit space
$M/G$) by the universal construction of Vinberg \cit!{32}: Define
the equivalence relation in $G \times C$ by
\begin{displaymath}
(x,g) \sim ( y,h)\;\;  \Longleftrightarrow \;\; x=y, \;\;  g^{-1}h \in G_x
\end{displaymath}
where $G_x =\langle s_{F_i}: x\in F_i \rangle$ is the subgroup
generated by all reflections with respect to walls containing $x$.
Then the quotient space
\begin{displaymath}
\mathcal U(G,C) = G \times C / \sim
\end{displaymath}
has the structure of a space of constant curvature such that the natural
action of $G$ on $\mathcal U(G,C)=M$ is isometric and $G$ is the reflection
group of $M$ with fundamental domain $C$.

More generally, if $G$ is a Coxeter group with a set
$S=\{s_1,\dots,s_l\}$ of standard generators, and relations
$s_i^2=1$, $(s_is_j)^{n_{i,j}}=1$, where $i,j=1,\dots,l$ and
$n_{ij}\in \mathbb N\cup \{\infty\}$, and if $C$ is a topological
space with closed subspaces $P_1,\dots,P_l$ 
(called panels), then
the Vinberg construction with $G_x =\langle s_{i}: x\in P_i \rangle$ 
gives a topological space $\mathcal U(G,C)$
with a continuous action of the group $G$ and orbit space $C$.
The topological $G$-space $\mathcal U(G,C)$ is called the {\it
universal space} of the Coxeter group $G$, and it satisfies the
following following universal property \cit!{32}: 
\begin{quote}
If $G$ acts in a topological space $X$ and if $\ph:C\to X$ is a
continuous map such that $s_i.\ph(x)=\ph(x)$ for $x\in P_i$ then
there exists a unique extension of $\ph$ to a $G$-equivariant
continuous map $\tilde\ph:\mathcal U(G,C)\to X$ such that
$\tilde\ph[1,x]=\ph(x)$ for any $x\in C$.
\end{quote}
Davis \cit!{12} found necessary and sufficient conditions that $\mathcal
U(G,C)$ is a topological manifold and $G$ is a topological reflection group
of $\mathcal U(G,C)$, i.e., any generator $s_i$ acts on $\mathcal U(G,C)$ as a
topological reflection (an involutive transformation whose fixed point set
$\mathcal U(G,C)^{s_i}$ separates $\mathcal U(G,C)$).
These conditions are that $C$ is a topological `nice' manifold
with corners and that each panel $P_i$ is a disjoint union of
walls such that for any $x\in C$ the subgroup $G_x=\langle s_i: x\in
P_i \rangle$ is finite. Conversely, let $G$ be a discrete group of
transformations of a topological manifold $M$ generated by
topological reflections, and let $C$ be its Weyl chamber (the
closure of a connected component of the set $M_{\text{reg}}=\{x\in
M:G_x=\{1\}\}$ of regular points). Let $s_1,\dots s_l$ be
reflections in $G$ such that $M^{s_i}\cap C$ contains a
codimension 1 component. Let $P_i$ be the union of all codimension
1 components $M^{s_i}\cap C$. Then $G$ is a Coxeter group with
standard generators $s_1,\dots,s_l$ and the $G$-manifold $M$ is
$G$-homeomorphic to the universal $G$-manifold $\mathcal U(G,C)$
defined by the panels $P_1,\dots, P_l$.

One of the aims of this paper is to describe the structure of the
Weyl chamber $C \cong M/G$ of a Riemannian manifold $M$ with a discrete
group $G$ generated by reflections, and to get a similar
description of such $G$-manifolds $M$ in terms of `abstract Riemannian
chambers' $C$, which are  Riemannian  manifolds with corners such that  
any two neighbouring walls $F_i, F_j$ satisfy the Coxeter property, i.e., 
the corresponding angle has constant value $\pi/n_{ij}$ along $F_i\cap F_j$.

In section \nmb!{2} we fix terminology and describe general properties of
reflections of a Riemannian manifold $M$ and of a discrete group $G$ generated 
by reflections. We discuss the relations between  a Dirichlet domain
$D$ of the group $G$ and its Weyl chamber $C$ which is defined as the closure of
a connected component of the set $M_{\text reg }$ of regular points of $G$.
We give an example when a Weyl chamber is larger than a Dirichlet domain. We
prove that for a simply connected manifold $M$, any reflection
$s$ is disecting, i.e. its fixed point set $M^s$ is a connected totally
geodesic hypersurface which decomposes $M$ into two parts.
We observe that a reflection group $G$ on a Riemannian manifold $M$
can be lifted canonically to a reflection group $\tilde G$,
which is an extension of $G$, on the universal covering $\tilde M$
of $M$. As an interesting example of Riemannian manifold with
a group generated by non disecting reflections, we consider
the maximal torus of the group $SU(n)$ for $n>2$ with the action
of the Weyl group.

Starting from section \nmb!{3}, we mostly consider a Riemannian
manifold $M$ with a reflection group $G$ generated by
{\it disecting} reflections. Such a $G$-manifold is called a {Coxeter
manifold}. Following M.\  Davis \cit!{12}, we derive from a lemma of
Bourbaki \cit!{4} that the group $G$ acts simply transitively on
the set of Weyl chambers $C$ of a Coxeter $G$-manifold $M$. This 
implies that Weyl chambers coincide with Dirichlet domains of regular
points and hence are homeomorphic to the orbit space, 
and that the reflection group $G$ is a Coxeter group with
reflections $s_i$ with respect to walls $F_i$ of $C$ as standard generators.
Moreover, the Weyl chamber $C$ has the structure of a Riemannian manifold
with corners and any two neighbouring walls $F_i, F_j$ of $C$ satisfy the Coxeter
property and yield a Coxeter relation $(s_is_j)^{n_{ij}}=1$. We prove that in the
simply connected case these relations generate all
relations of $G$. In the general case, we give a geometric description 
of the fundamental group $\pi_1(M)$.

In section \nmb!{4} we recall the notion and main the properties of a 
(smooth) manifold $M$ with corners and we define the concept of a Coxeter 
equipment of $M$. This is an order reversing mapping of the poset of
faces of $M$ into the poset of Coxeter subgroups of a given
Coxeter system $(G,S)$ (where  $S$ is the set of standard generators
of a Coxeter group $G$) which satisfies the Vinberg finiteness condition, see
\cit!{34}, \cit!{13}. We define a notion of {\it Riemannian chamber}
$C$ as a manifold with corners $C$ equipped with an
appropriate Riemannian metric such that walls $W_i$ of $C$
are totally geodesic and neighbouring walls satisfy the Coxeter property.
Any Riemannian chamber carries a universal Coxeter
equipment.

The Weyl chamber $C$ of a Coxeter $G$-manifold $M$ has
the natural structure of a Riemannian chamber with an
admissible (in some rigorous sense)
Coxeter equipment. Moreover, this equipment is universal
if and only if $\pi_1(M)= \pi_1(C).$  Conversely, if $C$
is a Riemmanian chamber with an admissible Coxeter $(G,S)$-equipment
then the universal space $M = \mathcal U(G,C)$ has the
structure of a Coxeter $G$-manifold with Weyl chamber
$C$. We prove also that if $C$ is a manifold with corners
and $s$ is a Coxeter equipment of $M$ then there exist a
Riemannian metric $\gamma$ such that $(M,\gamma)$ is a
Riemannian chamber and the equipment $s$ is admissible.
Hence any manifold with corners $C$ with a Coxeter
equipment determines a Coxeter $G$-manifold $M$, where
the metric of $M$ depends on the admissible metric on
$C$ and any Coxeter manifold can be obtained by this
construction.

In section \nmb!{5} and \nmb!{6} we discuss another approach for
reconstructing the Coxeter manifold from its Weyl
chamber $C$ which can be identified with the orbit
space $M/G$ based on the Thurston construction \cit!{31} of
the universal covering orbifold. We recall this
construction in section \nmb!{5} and we derive from the
main theorem of \cit!{23} that an orbifold structure
of a space $X$ can be reconstructed from the sheaf
$S_X$ of its smooth functions. In section \nmb!{6} we define
the notion of a Coxeter orbifold as an orbifold whose local
groups are finite linear Coxeter groups. An example of
Coxeter orbifold is the Weyl chamber $C$ of a Coxeter
manifold $M$.  We prove that any Coxeter orbifold is
such a Weyl chamber. More precisely, the universal
covering $M = \tilde C$ of a Coxeter orbifold $C$
admits a structure of (smooth) Coxeter $G$-manifold
such that $C$ is isomorphic to the Weyl chamber of
the isometry group $G$. In particular, this shows
that any Coxeter orbifold is {\it good} in the sense
of Thurston.

In the last section we described all Coxeter
equipments of an $n$-simplex $\Delta_n$. This gives a
classification of Coxeter orbifold structures on
$\Delta$ and a classification of Coxeter manifolds with 
orbit space $\Delta_n$ up to a diffeomorphism.

\section*{\totoc\nmb0{2}. Groups of
isometries generated by reflections}

\begin{proclaim}{\nmb.{2.1} Lemma}
Let $M$ be a connected complete Riemannian manifold, and let 
$G\subseteq  \on{Isom}(M)$  be a group of isometries.
Then $G$ is a discrete subgroup in the Lie group
$\on{Isom}(M)$ if and only if each orbit of $G$ in
$M$ is discrete. 
\end{proclaim}

We shall say that $G$ acts discretely on $M$.

\begin{demo}{Proof}
The pointwise-open topology on the Lie group
$\on{Isom}(M)$ of all isometries
coincides with the compact open topology. 

If $G$ is a discrete subgroup in
$\on{Isom}(M)$ then it is closed and acts
properly on $M$  so the action admits slices, and the orbit
$G.x$ through $x\in M$ is homeomorphic to $G/G_x$ where $G_x$
is the isotropy group of $x$. Thus each orbit is discrete.

Conversely, suppose that each orbit is discrete.
Since $G$ consists of isometries, each discrete orbit is closed. 
We consider the closure $\bar G$ of $G$ in $\on{Isom}(M)$. 
Since $G$-orbits are closed, $\bar G.x=G.x$ for each $x\in M$. 
The action of the closed group $\bar G$ of isometries is proper, 
so there exist slices. Let
$x_0$ be a regular point for the $\bar G$-action. Since $\bar
G.x_0$ is discrete, the slice $S_{x_0}$ through $x_0$ is open
in $M$, and the isotropy group $\bar G_{x_0}$ acts trivial on
$S_{x_0}$. Thus $\bar G_{x_0}$ acts trivial on $M$ and equals
$\{e\}$. Then $G.x_0=\bar G.x_0 \cong \bar G$, thus $\bar G=G$
and is discrete in $\on{Isom}(M)$.
\qed\end{demo}

\subsection*{\nmb.{2.2}. Dirichlet domains and central hypersurfaces}
Let $G\subset \on{Isom}(M)$ be a group which acts
isometrically and discretely 
on a connected complete Riemannian
manifold. 
Let $x_0$ be a regular point. The \idx{\it closed Dirichlet
domain} for this point is the set
\begin{displaymath}
D(x_0) := \{y\in M: d(y,x_0)\le d(y,g.x_0)
     \text{ for all }g\in G\},
\end{displaymath}
where $d$ is the geodesic distance on $M$. The open interior
$D(x_0)^o$ is called the \idx{\it open Dirichlet domain} for
the regular orbit $G.x_0$, and we can find a
\idx{\it fundamental domain} $F$ for the action of $G$ satisfying
$D(x_0)^o\subseteq F\subset D(x_0)$, i.e., a set
$F$ which meets each orbit in exactly one point, since
\begin{displaymath}
M=\bigcup_{g\in G}g.D(x_0).
\end{displaymath}
For any two different points $y_0,y_1\in M$ the \idx{\it central hypersurface}
is given by
\begin{displaymath}
H_{y_0,y_1} := \{y\in M: d(y,y_0) = d(y,y_1)\}.
\end{displaymath}
It disects $M$ in the sense that
$M\setminus H_{y_0,y_1}$ is the disjoint union of the two open sets
$\{x\in M: d(x,y_0)>d(x,,y_1)\}$ and $\{x\in M:d(x,y_0)<d(x,,y_1)\}$.
Note that if $M$ is a simply connected space of constant curvature
then $H_{y_0,y_1}$ is a totally geodesic submanifold, since it is the
fixed point set of a symmetry, but that in general 
$H_{y_0,y_1}$ is not a submanifold: On an elongated 2-torus it
can be a figure 8. 

\begin{proclaim}{Lemma}
For $x\in H_{y_0,y_1}$ let $c_0$ be a minimal geodesic from
$y_0$ to $x$. Then $c_0$ meets $H_{y_0,y_1}$ only at $x$.
\end{proclaim}

\begin{demo}{Proof}
Let $c_0(t_x)=x$ and suppose for contradiction that 
$c_0(t)\in H_{y_0,y_1}$ for $t<t_x$. Let $c_2$ be a minimal
geodesic from $x$ to $y_1$.
Then $t_x= d(y_0,c_0(t))+ d(c_0(t),x)= d(y_1,c_0(t))+
d(c_0(t),x)<d(y_1,x)=t_x$ unless $c_2$ equals the minimal
geodesic $s\mapsto c_1(t-s)$ and hence $y_0=y_1$, both a
contradiction.
\qed\end{demo}

\begin{proclaim}{\nmb.{2.3}. Lemma}
Let $D=D(x_0)$ be the closed Dirichlet domain of a regular point $x_0$ for a
discrete action of a group $G\subset \on{Isom}(M)$. Then
we have:
\begin{enumerate}
\item [(\nmb:{1})] If $g.D = D$ then $g=e$ in $G$.
\item [(\nmb:{2})] The open Dirichlet domain $D(x_0)^o$ is the
       connected component containing $x_0$ of
\begin{displaymath}
M\setminus \bigcup_{e\ne g\in G} H_{x_0,g.x_0}\subset M_{\text{reg}}.
\end{displaymath}
Here $M_{\text{reg}}$ denotes the set of all regular points, i.e.\ those points with trivial 
stabilizers. 
\item [(\nmb:{3})]
     $G$ acts simply transitively on the set $\{D(g.x_0):g\in G\}$ of all Dirichlet domains.
\end{enumerate}
\end{proclaim}

\begin{demo}{Proof} The isotropy group $G_{x_0}$ is trivial: See the proof of \nmb!{2.1}.

{\nmb|{1}}
If $g.D=D$ then $g.x_0\in D^o$, thus $d(g.x_0,x_0)\le d(g.x_0,h.x_0)$
for each $h\ne e$ in $G$. If $g\ne e$, putting $h=g$, we get
$g.x_0=x_0$, a contradiction.

{\nmb|{2}}
If $x\notin \bigcup_{e\ne g\in G}H_{x_0,g.x_0}$ then
$d(x,x_0)\ne d(x,g.x_0)$ for each $e\ne g\in G$. So if $g.x=x$ for
$g\ne e$ then
$d(x,x_0)=d(g.x,x_0)=d(x,g\i.x_0)\ne d(x,x_0)$, a contradiction. Thus the
isotropy group $G_x$ is trivial and $x$ is regular.
The connected component of $M\setminus \bigcup_{e\ne g\in G} H_{x_0,g.x_0}$
containing $x_0$ is the set of all $x\in M$ with $d(x,x_0)<d(x,g.x_0)$ for
all $e\ne g\in G$ which is $D^o$.

{\nmb|{3}}
Transitivity was seen in \nmb!{2.2} and simple transitivity follows from
{\nmb|{1}}.
\qed\end{demo}

\subsection*{\nmb.{2.4}. Walls of Dirichlet domains }
Let $G\subset \on{Isom}(M)$ be a discrete subgroup.
For a regular point $x_0\in M$ the set $F:=H_{x_0,g.x_0}\cap D(x_0)$
is called a \idx{\it wall} of the closed Dirichlet domain $D(x_0)$ if
it contains an open non-empty subset of $H_{x_0,g.x_0}$.

Two closed Dirichlet domains are called \idx{\it neighbors} if they
contain a common wall.

\begin{proclaim}{\nmb.{2.5}. Lemma of Poincar\'e}
Let $D=D(x_0)$ be a closed Dirichlet domain of a regular point $x_0$,
and let $g_1.D,g_2.D,\dotsc$ be all the neighbors of $D$. Then the
elements $g_1,g_2,\dotsc$ generate the group $G$.
\end{proclaim}

See fig.\ \fig!{4} for a Dirichlet domain with countably many
walls.

\begin{demo}{Proof}
{\it Claim.} For each $g\in G$ there exists a sequence
$e=h_0,h_1,\dots,h_n=g$ such that $D(h_i.x_0)$ and $D(h_{i+1}.x_0)$ are
neighbors for each $i$. We call this a Dirichlet neighbors chain
from $x_0$ to $g.x_0$.

The claim proves the lemma as follows. Since $D(h_1.x_0)$ is a neighbor of
$D=D(x_0)$ we have $D(h_1.x_0)=g_{i_1}.D$ for some $i_1$.
Then $g_{i_1}.g_{i_2}.D$ is the neighbor
$D(h_2.x_0)$ of $g_{i_1}.D$. Finally $g_{i_1}\dots g_{i_n}.D$ is the
neighbor $D(h_n.x_0)=D(g.x_0)$ of $g_{i_1}\dots
g_{i_{n-1}}.D=D(h_{n-1}.x_0)$. By \nmb!{2.3} we have
$g=g_{i_1}\dots g_{i_n}$.

We prove the claim by induction on $\{d_g:= d(x_0,g.x_0):g\in G\}$ which is
a locally finite set in $\mathbb R$ since the orbit $G.x_0$ is
discrete and closed in $M$.

Let $g\in G$ and assume that there exists a Dirichlet neighbors chain from
$x_0$ to $h.x_0$ whenever $d_h<d_g$. Applying $g_1$ we then also conclude
that there exists a Dirichlet neighbors chain from $g_1.x_0$ to $g_2.x_0$
whenever $d(g_1.x_0,g_2.x_0)<d_g$.
Consider a minimal geodesic $c$ from $x_0$ to $g.x_0$ of length $d_g$.

{\it Case \nmb:{1}.} Suppose that $c$ meets $\bigcup_{e\ne k\in G}H_{x_0,k.x_0}$
in $x=c(t_1)\in H_{x_0,k.x_0}$ at distance $t_1<\frac12 d_g$. See
fig.\ ~\fig!{1}.
\begin{figure}[h]
\centerline{
\psfrag{c}{$c$}
\psfrag{x}{$x$}
\psfrag{x0}{$x_0$}
\psfrag{g.x0}{$g.x_0$}
\psfrag{c1}{$c_1$}
\psfrag{c2}{$c_2$}
\psfrag{c3}{$c_3$}
\psfrag{k.x0}{$k.x_0$}
\psfrag{Hx0,k.x0}{$H_{x_0,k.x_0}$}
\includegraphics[width=8cm]{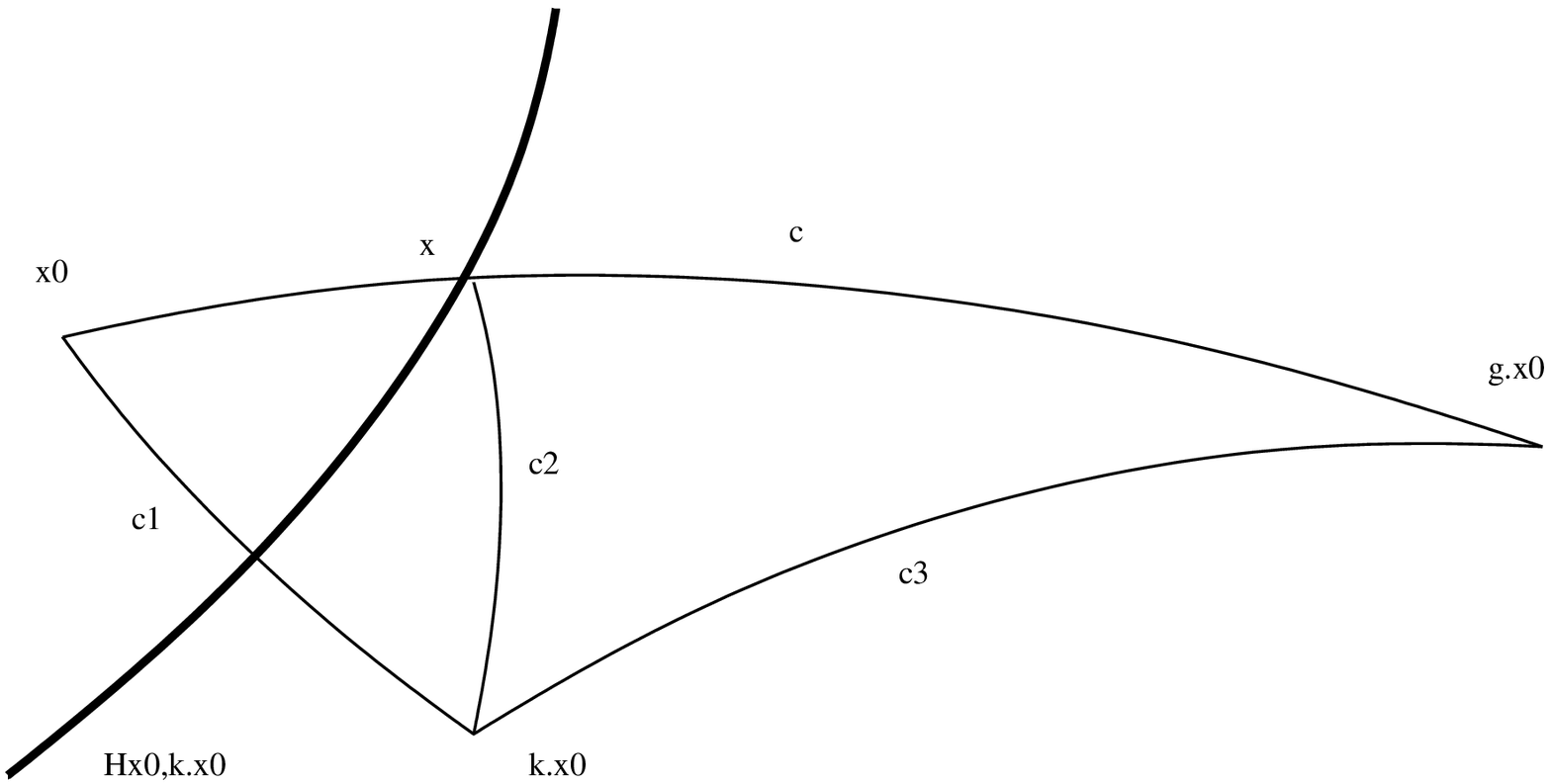}}
\caption{{\fig.{1}}}
\end{figure}\noindent
Then $k\ne g$. Consider a minimal geodesic $c_1$ from $x_0$ to $k.x_0$.
Then $c_1$ meets $H_{x_0,k.x_0}$ in $c_1(\frac12d_k)$. The minimal
geodesic $c_2$ from $x$ to $k.x_0$ has length $t_1$.
Thus $d_k\le 2t_1<d_g$ by the triangle inequality. By induction there
exists a Dirichlet neighbors chain from $x_0$ to $k.x_0$.
The minimal geodesic $c_3$ from $k.x_0$ to $g.x_0$ has lenght
$d_{k\i.g}< t_1 +d(x,g.x_0) = d_g$ since otherwise $k.x_0=x_0$. By
induction there exists a Dirichlet neighbors chain from $k.x_0$ to $g.x_0$.
So we get a Dirichlet neighbors chain from $x_0$ via $k.x_0$ to $g.x_0$, as
required.

{\it Case \nmb:{2}.} Suppose that $c$ meets $\bigcup_{e\ne k\in G}H_{x_0,k.x_0}$
for the first time at $x=c(\frac12d_g)\in H_{x_0,g.x_0}$ and that $x$ lies
in no other central hypersurface. Then there exists an open convex 
ball $U$ with center $x$ which meets only $H_{x_0,g.x_0}$ from 
$\bigcup_{e\ne k\in G}H_{x_0,k.x_0}$. 
By lemma \nmb!{2.2} the central hypersurface $H_{x_0,g.x_0}$ cuts $U$ in two connected
components $\{y\in U:d(x_0,y)\gtrless d(g.x_0,y)\}$ and is the boundary 
of both. One of them is in $D(x_0)^o$ and the other is in
the interior $D(h.x_0)^o$ of a neighbor of $D(x_0)$. 
So $y=c(\frac12d_g +\ep)\in D(h.x_0)^o$ for some $\ep>0$.
Then $d(y,h.x_0)<d(y,x_0)$ and
$d_{h\i g}=d(h.x_0,g.x_0)\le d(h.x_0,y)+d(y,g.x_0)< d_g$. By induction
there is a Dirichlet neighbors chain from $h.x_0$ to $g.x_0$, thus also
from $x_0$ to $g.x_0$.

{\it Case \nmb:{3}.} Suppose that $c$ meets $\bigcup_{e\ne k\in G}H_{x_0,k.x_0}$
for the first time at $x=c(\frac12d_g)\in H_{x_0,g.x_0}\cap H_{x_0,k.x_0}$
for $k\ne g$. We have $d(x,x_0)=d(x,k.x_0)$.
Consider the minimal geodesic $c_1$ from $x_0$ to $k.x_0$ which meets
$H_{x_0,k.x_0}$ in $c_1(\frac12 d_k)$. Then
$d_k=d(x_0,k.x_0)<d(x_0,x)+d(x,k.x_0)$ since otherwise the curve following
$c$ from $x_0$ to $x$ and then the minimal geodesic from $x$ to $k.x_0$
would be a minimal geodesic and could not have an angle $\ne 0$ at $x$
which implies that $k.x_0=g.x_0$. By induction there is Dirichlet neighbors
chain from $x_0$ to $k.x_0$. Moreover, $d_{k\i
g}=d(k.x_0,g.x_0)<d(k.x_0,x)+d(x,g.x_0)=d_g$ since otherwise the piecewise
minimal geodesic from $k.x_0$ via $x$ to $g.x_0$ would be a minimal
geodesic and thus $k.x_0=x_0$. By induction again there is a Dirichlet
neighbor chain from $k.x_0$ to $g.x_0$ which together with the first chain
gives a chain from $x_0$ to $g.x_0$, as required.
\qed\end{demo}
\subsection*{\nmb.{2.6}. Reflections }
Let $(M,\ga)$ be a connected complete Riemannian manifold. A
\idx{\it reflection} in $M$ is an isometry
$s\in \on{Isom}(M)$ such that for some fixed point $x_0$ of
$s$ the tangent mapping $T_{x_0}s$ is a reflection in the Euclidean
space $(T_{x_0}M,\ga_{x_0})$, with repect to a hyperplane:
For some vector
$0\ne X_{x_0}\in T_{x_0}M$ we have
$T_{x_0}s.X_{x_0}=-X_{x_0}$, whereas $T_{x_0}s|X_{x_0}^\bot = \Id$.

\begin{proclaim}{Lemma}
Let $s$ be a reflection on a complete connected Riemannian manifold
$M$. Then we have:
\begin{enumerate}
\item [(\nmb:{1})] Every connected component
       $N$ of the fixed point set $M^s$ is a totally geodesic
       submanifold, and for each $x\in N$ the tangent mapping $T_xs$
       equals the identity on $T_xN$ and $-\Id$ on $T_xN^\bot$.

\item [(\nmb:{2})] Every connected component $N$ of $M^s$ determines
       $s$ completely as follows: For $y\in M$ there exists $x\in N$ such
       that $d(y,x)=\on{dist}(y,N)$. Let
       $t\mapsto\exp(t.Y_x)$ be a minimal geodesic which reaches $y$
       at $t=1$. Then $s(y)=\exp(-Y_x)$.

\item [(\nmb:{3})] At least one connected
       component of $M^s$ is of codimension 1. Any such component is called
       a reflection hypersurface for $s$.

\item [(\nmb:{4})] For any $y\in M\setminus M^s$ we have
       $M^s\subseteq H_{y,s.y}$.
\end{enumerate}
\end{proclaim}

\begin{demo}{Proof}
{\nmb|{1}}
Let $x_0\in M^s$ be a point such that $T_{x_0}s$ is a Euclidean
reflection. Then $T_{x_0}s\o T_{x_0}s=\Id_{T_{x_0}M}$, thus $s$
is also an involution. Consequently $T_xs$ is an Euclidean involution for
each fixed point $x$, thus it is diagonalizable with eigenvalues $+1$
on the eigenspace $T_xN$ and eigenvalue $-1$ on the eigenspace
$T_xN^\bot$ where  $N$ is the connected component $N$ of
$M^s$ containing $x$.

{\nmb|{2}}
Note that $Y_x\in T_xN^\bot$ and that
$s(\exp(t.Y_x))=\exp(t.T_xs.Y_x)=\exp(-t.Y_x)$.

{\nmb|{3}}
The connected component of $M^s$ containing
$x_0$ is of codimension 1.

{\nmb|{4}} For $x\in M^s$ let $c$ be a minimal geodesic
from $x$ to $y$. Then $s\o c$ is a minimal geodesic from $s.x=x$ to
$s.y$. Thus $d(x,y)=d(x,s.y)$.
\qed\end{demo}

An example of a reflection $s$ which is generated by two
different reflection hypersurfaces $H,H'$: Let $M=S^1$, $H=\{1\}$,
$H'=\{-1\}$, and let $s$ be complex conjugation. Another
2-dimensional example with
three reflecting hypersurfaces is drawn in fig.\ \fig!{2}.
\begin{figure}[h]
\psfrag{1}{$1$}
\psfrag{-1}{$-1$}
\centerline{\includegraphics[width=8cm]{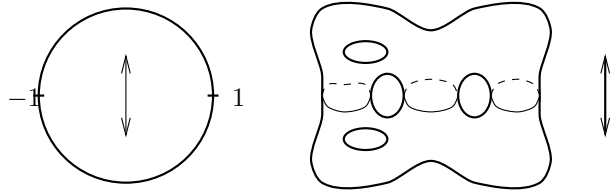}}
\caption{{\fig.{2}}
Reflections generated by different reflection hypersurfaces.}
\end{figure}\noindent
The 2-dimensional example in fig.\ \fig!{2} also shows that
two different reflection
hypersurfaces $H, H'$  for the same reflection $s$ need not be parallel, i.e.,
$\on{dist}(x,H')$ is not constant in $x\in H$.

In \cit!{30} one finds the following theorem: If an irreducible
Riemannian symmetric space $M$ of noncompact type admits an
involutory isometry whose fixed point set has codimension
one, then $M$ is a real hyperbolic space. This extends a result of
Iwahori \cit!{20} concerning irreducible Riemannian symmetric spaces
of compact
type; however, the proofs of these two results are substantially
different.

\subsection*{\nmb.{2.7}. Disecting reflections }
An isometry $s\in \on{Isom}(M)$ is called
\idx{\it disecting} if the complement of the 
fixed point set $M^s$ is not connected.

\begin{proclaim}{Lemma}
\begin{enumerate}
\item [(\nmb:{1})] A disecting isometry $s$ is a reflection.
\item [(\nmb:{2})] For a disecting reflection $s$ the fixed point set
       $M^s$ disects $M$ into exactly 2 pieces. The reflection $s$ permutes these
       two pieces.
\item [(\nmb:{3})] For a disecting reflection $s$ the
       fixed point set $M^s$ is a disjoint union of codimension 1
       submanifolds.
\item [(\nmb:{4})] For a disecting reflection $s$ and any
       $y\in M\setminus M^s$ we have $M^s=H_{y,s.y}$.
\end{enumerate}
\end{proclaim}

\begin{demo}{Proof}
{\nmb|{1}} Since $M\setminus M^s$ is not empty and
disconnected, the fixed point set $M^s$ which is a disjoint union of
closed totally geodesic submanifolds contains at least one
connected component of codimesion 1. For any $x$ in a codimension 1
component $H$ the tangent mapping $T_xs$ equals $\Id_{T_xH}$ on
$T_xH$ and is a nontrivial isometry on the 1-dimensional subspace
$T_xH^\bot$, thus equals multiplication by $-1$ there.

{\nmb|{2}}
Let $x_0\in M\setminus M^s$.
By \nmb!{2.6.4} we have $M^s\subseteq H_{x_0,s.x_0}$.
By \nmb!{2.3.2} for the group $\{\Id,s\}$, removing the set
$H_{x_0,s.x_0}$ decomposes $M$ into exactly two connected pieces.
Thus the subset $M^s\subseteq H_{x_0,s.x_0}$ cannot decompose it into
more than two pieces.

{\nmb|{3}} The union $M^s_1$ of all codimension 1 connected
components of $M^s$ also disects $M$ into two connected
components, since removing also the components of higher codimension
does not change connectedness any more.
Let $N$ be a connected component of codimension $\ge 2$
of $M^s$. Then $N$ is contained in one component of
$M\setminus M^s_1$ and $s$ thus has to map it into the other
component, by {\nmb|{2}}. Thus $N$ is empty.

{\nmb|{4}}
By \nmb!{2.6.4} we have $M^s\subseteq H_{y,s.y}$. Let
$z\in H_{y,s.y}$. If $z\notin M^s$ then $z$ and $s(z)$ lie in
different components of $M\setminus M^s$.
Let $c_1$ be a minimal geodesic from $y$ to $z$,
and let $c_2$ be a minimal geodesic (of the same length) from $s(y)$ to
$z$. Then the broken geodesic $c_1c_2\i$ from $y$ to $s(y)$ has to
meet $M^s$ in some point $x\in M^s$ since $y$ and $s(y)$ lie in
different components of $M\setminus M^s$.
If $x$ is an inner point on $c_2$, say, then the broken geodesic following
$s(c_2)$ from $y$ to $s(x)=x$ and then $c_2$ from $x$ to $z$ has the
same length as $c_2$ and hence $c_1$. It has an angle at $x$
(otherwise $z=s(z)$ and we are done), thus there is a geodesic from
$y$ to $z$ shorter than $c_1$, a contradiction.
\qed\end{demo}

\begin{proclaim}{\nmb.{2.8}. Theorem}
Let $M$ be a simply connected complete Riemannian
manifold.

Then
any reflection $\si$ on $M$ is disecting, and its fixed point set
$M^\si$ is a connected orientable totally geodesic closed
hypersurface.

\end{proclaim}

\begin{demo}{Proof} Let $x\in M\setminus M^\si$ and let $H$ be a connected
component of $M^\si$ of codimension 1.
Choose a minimal geodesic $c^+$ from $x$ to $H$. It hits $H$
orthogonally by minimality, and thus we may continue it by $c^-=\si.c^+$ to obtain a
geodesic $c_0$ from $x$ to $\si(x)$ which hits $H$ in exactly one
point.

Suppose that $M\setminus H$ is connected. Then there exists a smooth
curve $c_1$ from $x$ to $\si(x)$ in $M\setminus H$. Since $M$ is
simply connected, there exists a smooth homotopy
$h:[0,1]\x [0,1]\to M$ with $h(0,s)=x$, $h(1,s)=\si.x$,
$h(t,0)=c_0(t)$, and $h(t,1)=c_1(t)$. We can also assume that $h$ is
transversal to $H$. But then $h\i(H)$ is a closed 1-dimensional
submanifold in $[0,1]^2$ which hits the boundary exactly once in
$[0,1]\x \{0\}$ and never in $\{0\}\x [0,1]$ or $\{1\}\x [0,1]$. So the connected
component hitting once must hit again in $[0,1]\x \{1\}$. Thus $c_1$ hits
$H$, a contradiction.

Thus $M\setminus H$ is not connected, and $H$ cuts $M$ into two
components, $M^+$ and $M^-$. Moreover, $M^\si=H$ since $\si$
interchanges $M^+$ and $M^-$.
\qed\end{demo}

\subsection*{\nmb.{2.9}. Reflection groups and chambers}
Let $G\subset \on{Isom}(M)$ be a discrete subgroup of
isometries of a connected complete Riemannian manifold $M$
which is generated by all reflections contained in $G$.  We shall call any
such group $G$ a \idx{\it reflection group} of $M$.
By a \idx{\it (Weyl) chamber}  we mean the closure in
$M$ of a connected component of the
(open) complement of the union of all reflection hypersurfaces of all
reflections in $G$.
By an \idx{\it open (Weyl) chamber} we mean the open interior $C^o$ of a
Weyl chamber $C$.
For a chamber $C$ a \idx{\it wall} is a connected component of
$C\cap M^s$ for a reflection $s$ if it
contains a non-empty open subset of $M^s$ of codimension 1 in
$M$. Two walls  $F_i,F_j$ are called
neighbours is the intersection $F_i \cap F_j$ has a
connected component of codimension 2.

\begin{proclaim}{Lemma} For a chamber $C$ and any regular point $x_0\in C^o$
the Dirichlet domain $D(x_0)$ is contained in $C$.
Thus $C$ is a union of Dirichlet domains of the form $D(g.x_0)$,
for all $g\in N_G(C)$.

Moreover, $G$ acts transitively on the set of all chambers.
\end{proclaim}

\begin{demo}{Proof}
Since the set $M_{\text{reg}}$ of regular points is open and dense in $M$, 
we may choose a regular $x_0\in C$. We claim that $C\supseteq D(x_0)$.

For $x\in M$ consider a minimal geodesic $c$ from $x_0=c(0)$ to
$x=c(1)$. If $c$ hits a reflection hypersurface $H$ in $c(t)$ for
$t<1$, we may consider the minimal geodesic from $x_0$ to $c(t)$,
followed by the minimal geodesic from $c(t)$ to $s_H(x)$, which is a
broken geodesic from $x_0$ to $s_H(x)$ of length
$d(x_0,c(t))+d(c(t),s_H(x)) = d(x_0,c(t))+d(c(t),x) = d(x_0,x)$.
Since it has a proper angle at $c(t)$, a minimal geodesic from
$x_0$ to $s_H(x)$ has length $d(x_0,s_H(x))<d(x_0,x)$.
Thus we see: {\it Whenever the minimal geodesic from $x_0$ to a point
$x$ hits a reflection hypersurface $H$ in an intermediate point,
$d(x_0,s_H(x))<d(x_0,x)$.}

Now let $y\in D(x_0)$. Then
$d(x_0,y)=\on{dist}(x_0,G.y)$.
By the statement above, any minimal geodesic from $x_0$ to $y$ can
hit a reflection hypersurface at most in $y$. Thus $y\in C$.

That $G$ acts transitively on the set of all chambers follows from \nmb!{2.3}.
\qed\end{demo}

\subsection*{\nmb.{2.10}. Examples of non-disecting reflections }
We consider the real projective plane $\mathbb R\mathbb P^2$, with 
the metric induced from $S^2$, and a
reflection $s=s_H$ on one line $H$ in it.
Look at fig.\ \fig!{3} where the line at infinity $L$ is
chosen orthogonal to $H$ so that $L$ is invariant under $s$.
On the line $L$ one still has to
identify antipodically. The fixed point set of $s$ consists of $H$
and the single point $x$ on $L$ farthest from $H$.
There is only one chamber
$C=\mathbb R\mathbb P^2\setminus (H\cup\{x\})$ which is a punctured disk and is dense in
$\mathbb R\mathbb P^2$. But there are two Dirichlet domains
$D(z)$ and $D(sz)$ depending on $z$ which for $z$ on the line through $x$
othogonal to $L$ meet in $H\cup L$.
\begin{figure}[h]
\psfrag{L}{$L$}
\psfrag{H}{$H$}
\psfrag{x}{$x$}
\psfrag{s}{$s$}
\centerline{\includegraphics[width=4cm]{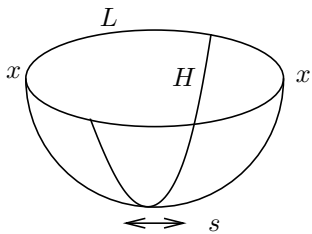}}
\caption{{\fig.{3}}
A reflection on $\mathbb RP^2$.}
\end{figure}\noindent
As another example,
consider $M=SO(3)=\mathbb RP^3$ with the biinvariant metric. Then $g\mapsto g\i$
is a non-disecting reflection whose fixed point set is the disjoint
union of $\{e\}$ and some $\mathbb RP^2$. This reflection generates a
Coxeter group.

\begin{proclaim}{\nmb.{2.11}. Theorem} Let $M$ be a complete Riemannian
manifold and let $G\subseteq \on{Isom}(M)$ be a discrete
group of isometries which is generated by all its reflections. Let
$C$ be a Weyl chamber in $M$ for $G$. Let $F_1,F_2,\dotsc$ be the
walls of $C$ and let $s_i $ be the reflection with respect to the
wall $F_i$.

Then the reflections $s_1,s_2,\dotsc$ generate
$G$, and they satisfy the following relations:
\begin{enumerate}
\item [(\nmb:{1})] $(s_i)^2=1$
\item [(\nmb:{2})] If  two walls $F_i, F_j$ are neighbors 
  then $(s_is_j)^{n_{ij}}=1$ for some natural number $n_{ij}$.
\end{enumerate}
\end{proclaim}

See fig.\ \fig!{4} for a Weyl chamber in the Poincar\'e
upper halfplane with infinitely many walls.

\begin{figure}[h]
\centerline{
\includegraphics[width=8cm]{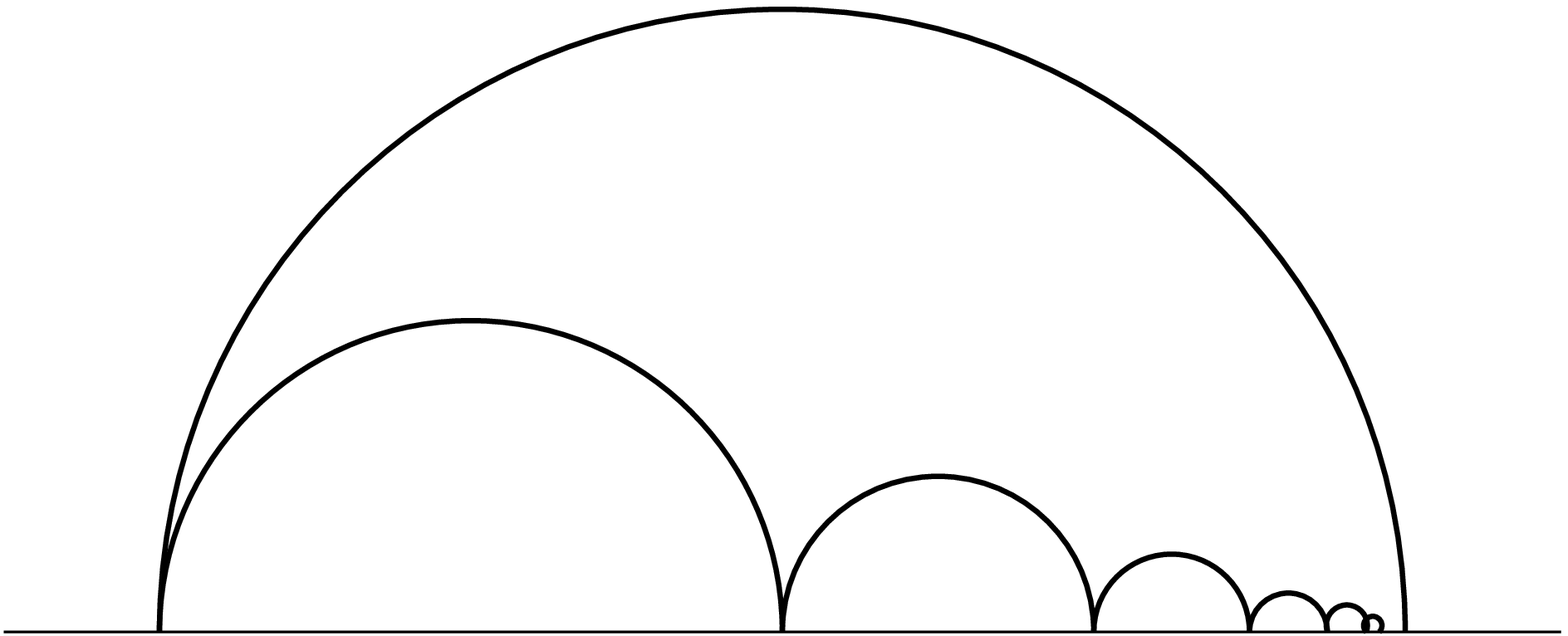}}
\caption{{\fig.{4}}} A Weyl chamber in the Poincar\'e
upper halfplane with infinitely many walls.
\end{figure}\noindent

\begin{demo}{Proof}
We prove that the reflections
$s_1,s_2,\dotsc$ generate $G$.
Let $C'$ be any other Weyl chamber in $M$. Then we choose a smooth curve
$c:[0,1]\to M$ from a regular point $x_0\in C$ to a regular point
$x'\in C'$ which changes Weyl chambers only transversally through
open interiors of walls. First the curve passes from $C$ through the
interior of a wall $F_{i_1}$ to a neighbor $s_{i_1} (C)$, and then
through a wall $F$ of this chamber to the next. For the reflection $s_F$ in
$F$ we have $s_{F}=s_{i_1}. s_{i_2} . s_{i_1 }$
for some wall $F_{i_2}$ of $C$.
If we now follow the curve $c$ through all interiors of walls we see
that $C'$ is of the form $C'=g(C)$ for $g$ in the subgroup generated
by $s_1 ,\dots,s_l$.
Any reflection in $G$ is of the form $s_{F}$ for some wall of some
chamber $C'$. But then the argument above shows that
$s_F=g. s_{i_k}. g\i$, so $G$ is generated
by $s_1 ,s_2,\dotsc$ as claimed.

Relations {\nmb|{1}} and {\nmb|{2}}
follow, since if $x$ is an interior point of the face $f=
F_i \cap F_j$
(i.e., there are no other walls through $x$) then the stabilizer
$G_x$ is faithfully and orthogonally represented in the two-dimensional
space $T_x(f)^{\perp}$, and any finite subgroup of $O(2)$ which contains
a reflection is a dihedral group.
\qed\end{demo}

\subsection*{\nmb.{2.12}. Remark} In the setting of theorem \nmb!{2.11}
there might be more relations than specified in \nmb!{2.11.1} and
\nmb!{2.11.2}, see fig.\  \fig!{5}.
\begin{figure}[h]
\psfrag{s1}{$s_1$}
\psfrag{s2}{$s_2$}
\psfrag{s3}{$s_3$}
\psfrag{s4}{$s_4$}
\psfrag{C}{$C$}
\centerline{\includegraphics[width=8cm]{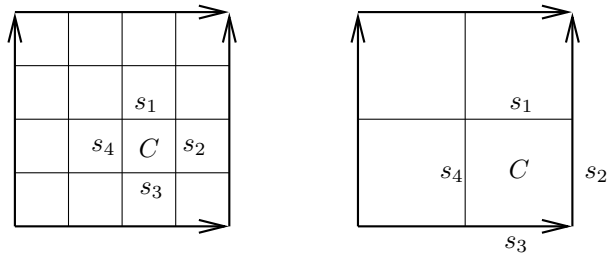}}
\caption{{\fig.{5}}
2-tori with $\mathbb Z_2^4$ and $\mathbb Z_2^2$ as reflection
groups.}
\end{figure}\noindent
The left part of fig.\  \fig!{5} is a  flat 2-torus
with a chamber $C$ specified, with reflections $s_1$, $s_2$, $s_3$, $s_4$
and angular relations $(s_is_{i+1})^2=1$ for $i=1,\dots,4 \mod(4)$ as
decribed in \nmb!{2.11.2}.
But moreover the relations $(s_2s_4)^2=1$ and $(s_1s_3)^2=1$ hold
which are not described by \nmb!{2.11.2}.

In right hand part of fig.\  \fig!{5} we even have $s_1=s_3$ and $s_2=s_4$.

\subsection*{\nmb.{2.13}. Lifting reflection groups to the universal covering}
Let $\pi:\tilde M\to M$ be the universal covering of a Riemannian manifold
$M$ with a reflection group $G$, and let 
$\pi_1(M)=\Ga\subset \on{Isom}(\tilde M)$
be the group of deck transformations of $\pi$. Any isometry of $M$ can be lifted
to an isometry of $\tilde M$. A lift $\tilde s$ of a reflection $s$ in $G$
is a reflection on $\tilde M$ if
and only if it has a fixed point $\tilde x\in \tilde M$ with $\pi(\tilde
x)$ in a reflection hypersurface of $s$ in $M$. The group $\tilde G$
generated by all reflections which are lifts of reflections in $G$, is a
reflection group in $\tilde M$ which is normalized by $\Ga$ in $\on{Isom}(\tilde
M)$. Then
$\tilde G \Ga$ is the group of all lifts of transformations in $G$, and
$G=(\tilde G\Ga)/\Ga=\tilde G/(\tilde G\cap \Ga)$.
If $\tilde C$ is a chamber for $\tilde G$ in $\tilde M$ then $\pi(\tilde
C)$ is a chamber for $G$ in $M$, since the union of all
reflections hypersurfaces of $\tilde G$ equals the inverse image under $\pi$ of 
the union of all reflection hypersurfaces of $G$. 

Let $s$ be a reflection in $G$, and let $\tilde s$ be
a reflection covering $s$ in $\tilde G$. According to \nmb!{2.8} 
each reflection $\tilde s$ in $\tilde G$
is disecting, $\tilde M^{\tilde s}$ is one reflection
hypersurface, and $\tilde M\setminus \tilde M^{\tilde s}$
consists of exactly two connected components 
$\tilde M^{\tilde s}_+$ and 
$\tilde M^{\tilde s}_-$. 

{\it If $G$ is generated by disecting reflections 
then $G$ acts simply transitively on the set of all chambers},
see \nmb!{3.5} below.  
The converse is not true, even if $G$ is a Coxeter group, 
see fig.\ \fig!{7} in \nmb!{2.15}.

{\it Suppose that one (equivalently any) chamber is simply
connected. Then $G$ acts simply transitively on the set of all
chambers if and only if $\Ga\subseteq \tilde G$.} To see this,
note that the universal cover $\pi:\tilde M\to M$ restricts to a
diffeomorphism for each chamber $\tilde C$ in $\tilde M$ onto a
chamber $C=\pi(\tilde C)$ in $M$. 
If $\Ga$ contains a nontrivial deck transformation $\ga$, then
for a chamber $\tilde C$ covering $C$ the set $\ga(\tilde C)$
is another chamber covering $C$. By \nmb!{2.10} and \nmb!{3.5}
there exists a unique $\tilde g\in \tilde G$ with $\tilde
g(\tilde C)=\ga(\tilde C)$. But then $\tilde g=\ga$ 
if and only if $\tilde g$ covers $\on{Id}_M$ in $G$.

\begin{proclaim}{\nmb.{2.14}. Proposition}
Let $G$ be a reflection group on a simply connected 
complete Riemannian manifold $M$.
Then each chamber $C$ is simply
connected.
\end{proclaim}

\begin{demo}{Proof}
Suppose for contradiction that some chamber $C$ is not simply
connected:
Let $c:[0,1]\to C$ be a closed smooth curve through a
regular point $x_0\in C$ which is not contractible to the constant
curve through $x_0$ in $C$ with fixed ends at $x_0$.

Since $M$ is simply connected there exists a smooth homotopy
$h:[0,1]\x [0,1]\to M$ with $h(0,t)=c(t)$, $h(s,0)=x_0$,
$h(s,1)=x_0$, and $h(1,t)=x_0$.
We may assume that $h$ is transversal to each reflection hypersurface
and to each intersection of such hypersurfaces, since these form a
locally finite family by the discreteness of $G$. Thus for each
intersection hypersurface $H_i$ the set $h\i(H_i)$ is a 1 dimensional
embedded submanifold of $[0,1]^2$ which does not meet the boundary,
so it is a disjoint set of embedded circles in $C$ which may
touch only the bottom boundary $\{0\}\x[0,1]$.
Moreover, the sets $h\i(H_i)$ are all pairwise transversal
1-dimensional submanifolds in $(0,1)^2$, or empty,
since this is the case for the (geodesically closed) $H_i$ in $M$.
Fig.\ \fig!{6} is an illustration. 
See \cit!{25}, section 6, for transversality theorems on manifolds
with corners. 
\begin{figure}[h]
\psfrag{s}{$s$}
\psfrag{t}{$t$}
\centerline{\includegraphics[width=4cm]{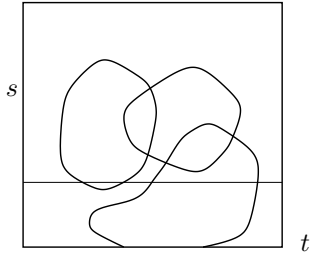}}
\caption{{\fig.{6}}
The sets $h\i(H_i)$ in $[0,1]^2$.}
\end{figure}\noindent
Now $h_0=c$ is completely contained in $\bar c$ and we consider the curve
$h_s=h|(\{s\}\x [0,1])$ for $s$ moving from 0 to 1. So we move
$\{s\}\x[0,1]$ upwards inside $[0,1]^2$.
If this line hits $h\i(H_i)$
we start reflecting back into $C$ the point $h_s(t)$ for those
$t$ which lie inside $h\i(H_i)$. If we meet another $h\i(H_j)$ we add
the reflection $s_{H_j}$ at the right, etc. Since the different
$h\i(H_i)$ are transversal to each other this is welldefined, in
particular at $s=0$, where the sets $(\{0\}\x [0,1])\cap h\i(H_i)$
are disjoint by transversality.
This proceedure transforms the smooth homotopy $h:[0,1]^2\to M$ to a
continuous homotopy $\bar h:[0,1]^2\to C$ which contracts $c$ to
$x_0$. Thus $C$ is simply connected which contradicts our
assumption.
\qed\end{demo}

\subsection*{\nmb.{2.15}.
Maximal torus of a compact Lie group as manifold with reflections
}
Let $G$ be a semisimple compact Lie group with Lie algebra $\mathfrak g_0$ and
let $T$ be a maximal torus in $G$. The Lie subalgebra $\mathfrak t_0$ to $T$ is
then a Cartan subalgebra. Let $\De\subset \mathfrak t^*$ be the set of roots
where $\mathfrak t= \mathfrak t_0\otimes \mathbb C$ is the complexification of 
$\mathfrak t_0$ and where $\mathfrak t^*$ is the dual space of $\mathfrak t$. 
Each root is
purely imaginary on $\mathfrak t_0$. We have the following inclusion of
lattices in $\mathfrak t^*$:
\begin{displaymath}
\mathbb Z\De \subseteq \La_{\text{anal}} \subseteq \La_{\text{alg}},
  \qquad\text{ where}
\end{displaymath}
\begin{enumerate}
\item [$\mathbb Z\De$] is the root lattice, generated by $\De$,
\item [$\La_{\text{anal}}$] is the lattice of analytically integral forms
  $\la\in L(\mathfrak t_0,i\mathbb R)$; 
  they are characterized by the following property:
  whenever $H\in \mathfrak t_0$ satisfies $\exp(H)=1$ then $\la(H)\in 2\pi
  i\mathbb Z$; equivalently: there exists a multiplicative character
  $\xi_\la:T\to S^1$ such that $e^{\la(H)}=\xi_\la(\exp(H))$ for all
  $H\in\mathfrak t_0$.
\item [$\La_{\text{alg}}$]
  is the weight lattice consisting of all 
  $\la\in L(\mathfrak t_0,i\mathbb R)$ such that 
  $2\langle\al,\la\rangle/|\al|^2\in \mathbb Z$ for all roots $\al\in\De$.
\end{enumerate}
Now $\exp:\mathfrak t_0\to T$ induces an isomorphism
$\mathfrak t_0/\La_{\text{anal}}^*=T$, where $\La_{\text{anal}}^*$ is the
dual lattice 
$\{X\in \mathfrak t_0:\la(X)\in \mathbb Z\text{ for all }\la\in\La_{\text{anal}}\}$.
Recall that $G$ has trivial center if and only if
$\La_{\text{anal}}=\mathbb Z\De$, that $G$ is simply connected if and only
if $\La_{\text{anal}}=\La_{\text{alg}}$, that in general
$\La_{\text{anal}}/\mathbb Z\De$ is the center of $G$, and that the order of
$\La_{\text{alg}}/\mathbb Z\De$ equals the determinant of the Cartan matrix of
$\g$.
\begin{figure}[h]
\psfrag{1}{1}
\psfrag{2}{2}
\psfrag{3}{3}
\psfrag{4}{4}
\psfrag{5}{5}
\psfrag{6}{6}
\centerline{\includegraphics[width=8cm]{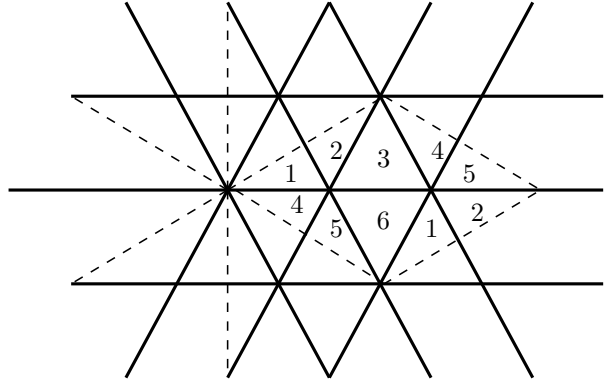}}
\caption{{\fig.{7}} 
$\mathfrak t_0$ as universal covering of $T$ for $A_2=SU(3)$.
}
\end{figure}\noindent
The reflections on $T$ are induced by the reflections in the Weyl group in
$\mathfrak t_0$; to visualize it we consider the reflections hyperplanes and the
lattice $\La_{\text{anal}}^*$ which consists of vectors orthogonal to the
reflection hyperplanes. Then we consider a standard fundamental domain of
the additive action of $\La_{\text{anal}}^*$.
We see that for $A_1^k=SU(2)^k$ all reflections in $T$ are disecting,
but that for semisimple nonabelian $G$ we always get nondisecting
reflections. 

See fig.\ \fig!{7} for an example: It shows for 
$A_2=SU(3)$ the Cartan algebra $\mathfrak t_0$ as the universal covering of $T$ 
with the reflection hyperplanes (bold) for $W\rtimes \La^*_{\text{anal}}$, 
the lattice $\La^*_{\text{anal}}$, and the fundamental domain (dashed).
The reflections on $T$ are not disecting, 
and the reflection group acts freely on the set of 
chambers in $T$, which are numbered.

\section*{\totoc\nmb0{3}. Coxeter Riemannian manifolds }

\subsection*{\nmb.{3.1}. Coxeter groups }
\cit!{4}
Recall that a \idx{\it Coxeter group} is a group $G$ which is a quotient
of a free group $G(S)$ with a set $S$ of
generators by the subgroup generated by the relations $s^2=1$ and
$(ss')^{n_{s,s'}}=1$ for all $s,s'\in S$, where
$n_{s,s'}\in\{1,2,\dots,\infty\}$ indicates the order of $ss'$ in $G$.

The set $S$ is called a set of \idx{\it standard generators} of $G$, and
$(G,S)$ is called a \idx{\it Coxeter system} for $G$.
Any subset $S'\subset S$ generates a subgroup $G(S')\subset G$ such that
$(G(S'),S')$ is again a Coxeter system. $G(S')$ is called a {\it Coxeter
subgroup}. The set of all Coxeter subgroups is a partially ordered set with
respect to inclusion.
A Coxeter system is described by a Coxeter diagram with vertices
corresponding to the elements of $S$, where $s$ and $s'$ are connected by
$n_{s,s'}-2$ edges if $(ss')^{n_{ss'}}=1$ and $1<n_{ss'}<\infty$.
The Coxeter diagram of a Coxeter subgroup $(G(S'),S')$ for $S'\subset S$ is
obtained from the Coxeter diagram of $(G,S)$ by deleting all vertices in
$S\setminus S'$ and all edges leading to such vertices.

The length $\ell(g)$ of an element $g\in G$ is the minimum number $l$
such that $g=s_{i_1}\dots s_{i_l}$ for $s_{i_k}\in S$. It satisfies
$\ell(gg')\le \ell(g)+\ell(g')$, $\ell(g\i)=\ell(g)$, and
$|\ell(g')-\ell(g)|\le \ell(g'g\i)$.

In a Coxeter group $(G,S)$ let $P^+_s:=\{g\in G: \ell(sg)>\ell(g)\}$
and $P_s^-:=sP_s^+$. Then we have \cit!{4},~iv,~1,~7:
\begin{enumerate}
\item[(\nmb:{1})] $\bigcap_{s\in S} P_s^+=\{e\}$.
\item[(\nmb:{2})] $G= P_s^+\sqcup P_s^-$ (disjoint union) for each $s\in S$.
\item[(\nmb:{3})] Let $s,s'\in S$ and $g\in G$. If $g\in P_s^+$ and
       $gs'\notin P_s^+$ then $s=gs'g\i$.
\end{enumerate}

Conversely, let $G$ be a group with a generating set $S$ of
idempotents. Let $(P_s)_{s\in S}$ be a family of subset of $G$ which
satisfies
\begin{enumerate}
\item[(\nmb:{4})] $e\in P_s$ for all $s\in S$.
\item[(\nmb:{5})] $P_s\cap sP_s = \emptyset$ for all $s\in S$.
\item[(\nmb:{6})] Let $s,s'\in S$ and $g\in G$. If $g\in P_s$ and
       $gs'\notin P_s$ then $s=gs'g\i$.
\end{enumerate}
Then $(G,S)$ is a Coxeter system and $P_s=P_s^+$.

\subsection*{\nmb.{3.2}. Riemannian Coxeter manifold }
Let $G\subset \on{Isom}(M)$ be a discrete subgroup of
isometries of a complete  Riemannian manifold $M$
which is generated by disecting reflections. Then $(M,G)$ is called a
{\it Riemannian Coxeter manifold}.

\subsection*{\nmb.{3.3}. 
Coxeter manifolds of constant curvature}
We recall some classical results.

Let $(G,S)$ be a Coxeter system such that $G$ is a finite group and let
$S=\{s_1,\dots,s_n\}$.
Then there exists a unique orthogonal representation of $G$ as a linear reflection
group on an Euclidean space $\mathbb R^n$ such that the $s_i$ are
reflections. The Weyl chamber associated to $S$
is a simplicial cone with walls $F_1,\dots, F_n$ such that $s_i$ is the
reflection in $F_i$.
Then the angle $\al_{ij}$ between $F_i$ and $F_j$ is given by
$\al_{ij}=\pi/n_{ij}$ where $(s_is_j)^{n_{ij}}=1$ and $n_{ij}$ is minimal.
In the following table we give the list of all finite Coxeter systems which are
irreducible in the sense that they are not a direct product of two
(commuting) Coxeter subsystems.
\allowdisplaybreaks
\begin{alignat*}3
&A_l &\qquad&l\ge1 &\qquad &\dynkin\root{s_1}\link{}\root{s_2}\link{}\dots
  \link{}\root{s_l}\enddynkin\\
&B_l &\qquad&l\ge2 &\qquad &\dynkin\root{s_1}\link{}\root{s_2}\link{}\dots
  \link{}\root{s_{l-1}}\llink{}\root{s_l} \enddynkin\\
&D_l &\qquad&l\ge4 &\qquad
  &\dynkin\root{s_1}\link{}\root{s_2}\link{}\dots\link{}\rroot{s_{l-2}}
  \rootupright{s_{l-1}}\rootdownright{s_l} \enddynkin\\
&E_6  &\qquad& &\qquad
  &\dynkin\root{s_1}\link{}\root{s_2}\link{}\root{s_3}
  \rootdown{s_4}\link{}\root{s_5}\link{}\root{s_6}\enddynkin\\
&E_7  &\qquad& &\qquad
  &\dynkin\root{s_1}\link{}\root{s_2}\link{}\root{s_3}\link{}\root{s_4}
  \rootdown{s_5}\link{}\root{s_6}\link{}\root{s_7}\enddynkin\\
&E_8  &\qquad& &\qquad
  &\dynkin\root{s_1}\link{}\root{s_2}\link{}\root{s_3}
  \link{}\root{s_4}\link{}\root {s_5}
  \rootdown{s_6}\link{}\root{s_7}\link{}\root{s_8}\enddynkin\\
&F_4 &\qquad& &\qquad &\dynkin\root{s_1}\link{}\root{s_2}
  \llink{}\root{s_3}\link{}\root{s_4} \enddynkin\\
&\on{Di}(k+2) &\qquad&k\ge 4 &\qquad &\dynkin\root{s_1}\link{^k}\root{s_2}
  \enddynkin\qquad\qquad\qquad\qquad\qquad  G_2=\on{Di}(6)\\
&H_3 &\qquad&m\ge 3 &\qquad &\dynkin\root{s_1}\link{^3}\root{s_2}
  \link{}\root{s_3}\enddynkin\\
&H_4 &\qquad&m\ge 3 &\qquad &\dynkin\root{s_1}\link{^3}\root{s_2}
  \link{}\root{s_3}\link{}\root{s_4}
  \enddynkin\\
\end{alignat*}
If the Coxeter group has no dihedral group $\on{Di}(k+2)$ as direct factor,
then the angle between two walls may only take the values $\al=\pi/n$
for $n=2,3,4,5,6$.

Conversely any simplicial cone with walls $F_1,\dots,F_n$ having angles
$\al_{ij}=\pi/n_{ij}$ between $F_i$ and $F_j$ where $n_{i,j}\in \mathbb N$, is
the Weyl chamber of a uniquely given Coxeter system with finite Coxeter
group, by \cit!{32}, theorem 1.
The Coxeter diagram of $(G,S)$ contains also all information about the Weyl
chamber. The angle between the walls $F_i$ and $F_j$ is
$\al_{ij}=\pi/n_{i,j}$ where $n_{ij}-2$ is the number of edges connecting
the vertices $s_i$ and $s_j$.

If $g\in G$ preserves a codimension $k$ face (an intersection of
$k$ walls) $F=F_{i_1}\cap\dots\cap F_{i_k}$ which does not contain a line
through 0,
then it it preserves it pointwise. Namely, $g$ has a fixed point $x$ in the
interior of $F$ since $F$ is convex. By
the lemma of Chevalley, $g$ is contained in the Coxeter subgroup generated
by all reflections $s_i$ fixing $x$ which correspond to all walls
through $x$. Since $x$ is an inner point of $F$, these walls also contain $F$.
Thus $g$ fixes $F$ pointwise.

The angle in $F_{i_1}\cap\dots\cap F_{i_k}$ between
$F_{i_1}\cap\dots\cap F_{i_k}\cap F_{i_{k+1}}$ and $F_{i_1}\cap\dots\cap
F_{i_k}\cap F_{i_{k+2}}$ is in general not of the form $\pi/n$;
nevertheless it is uniquely determined by the Coxeter system.

\subsection*{\nmb.{3.4}. Example of Coxeter manifolds of non-constant
curvature}
Let $G$ be a linear reflection group on $\mathbb R^n$.

(\nmb:{1}) 
Let $S$ be the unit sphere of $\mathbb R^n$. Then $G$ acts on $S$ and is
generated by reflections. Choose a chamber $C$ in $\mathbb R^n$ and a
$(n-1)$-ball $B$ in $C\cap S$. By surgery one may glue any compact
$(n-1)$-dimensional manifold $M$ to $\partial B$ and do this in each
chamber via the transformations of the group $G$. Obviously one can
also put a $G$-invariant Riemannian metric on the resulting
manifold, which then has complicated topology but carries a finite
subgroup of the group of isometries which is generated by
disecting reflections.

(\nmb:{2})
Choose a chamber $C$ in $\mathbb R^n$ and within $C$ a regular point.
Connect this point by a smooth curve to some point in each interior
of each wall of $C$. Distributing this by $G$ into all chambers of
$\mathbb R^n$ yields a graph on which $G$ acts. Now replace each point in the
walls by a $S^1$ which lies completely in the interior of the wall,
and replace the piece of the graph in the chamber $C$ by a smooth
compact surface which all the $S^1$'s as boundary components, meeting
the walls orthogonally. Distribute this to all chambers by the
$G$-action and obtain a
smooth compact surface with induced Riemannian metric on which $G$
acts as a group of isometries generated by reflections.

\begin{proclaim}{\nmb.{3.5}. Theorem}
Let $(M,G)$ be a Riemannian Coxeter manifold.
Then $G$ is a Coxeter group and $(G,S)$ is a Coxeter system for $G$,
where $S$ is the set of reflections with respect to the walls
of $C$.
Moreover, $G$ acts simply transitively on the set of chambers.
\end{proclaim}

\begin{demo}{Proof} We follow arguments from \cit!{12}. 
Let $Q$ be a chamber. For a reflection $s$ with respect to a wall $F$
of $Q$ we set
\begin{displaymath}
P_s := \{g\in G: gQ\subset M^s_+\}
\end{displaymath}
where $M^s_+ $ is the connected componet of $M \setminus M^s$ which contains
$Q$.

\begin{proclaim}{Lemma}
$P_s = P_s^+ = : \{ g \in G : \ell(sg) > \ell (g) \}.$
\end{proclaim}
\begin{demo}{Proof} It is sufficent to check the properties
\nmb!{3.1.4}, \nmb!{3.1.5}, and \nmb!{3.1.6}.
The first two properies are obvious. We check \nmb!{3.1.6}.
Let $s,s'$ be reflections with
respect to walls $F,F'$ of the chamber $Q$ and $g \in P_s$ but $gs' \notin P_s$.
The chambers $Q,s'Q$ have a common wall $W$  and the chambers $gQ$, $gs'Q$
have a common wall $gW$.
Since they are on different sides of the hypersurface  $M^s$, the wall
$gW$ belongs to $M^s$, see Fig. \fig!{8}.
\begin{figure}[h]
\psfrag{s'Q}{$s'Q$}
\psfrag{W}{$W$}
\psfrag{Q}{$Q$}
\psfrag{s}{$s$}
\psfrag{s'}{$s'$}
\psfrag{gQ}{$gQ$}
\psfrag{gW}{$gW$}
\psfrag{gs'Q}{$gs'Q$}
\psfrag{Ms}{$M^s$}
\psfrag{Ms'}{$M^{s'}$}
\psfrag{Ms+}{$M^s_+$}
\psfrag{Ms-=sMs+}{$M^s_-=sM^s_+$}
\centerline{\includegraphics[width=10cm]{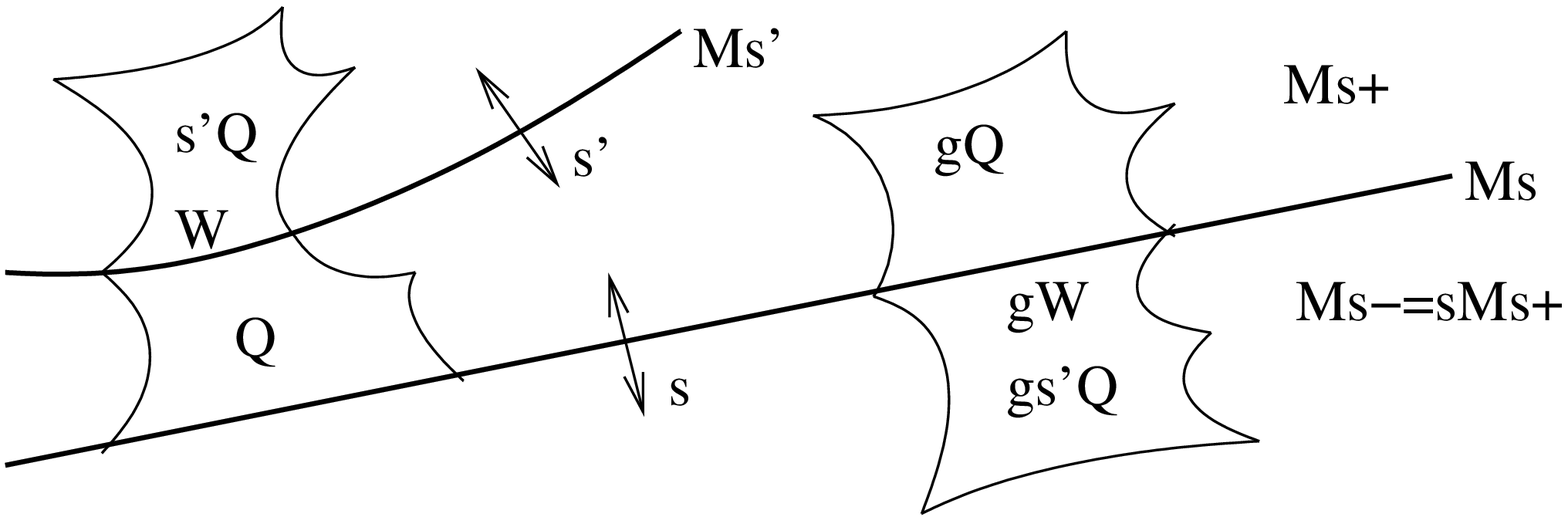}}
\caption{{\fig.{8}}
}
\end{figure}\noindent
Then  $s(gQ) = gs'Q $  and
$s'$, $g\i sg$ are two reflections which map $Q$ to $s'Q$.
Moreover $W\subset M^{s'}$, and $gW\subset M^s$ implies $W\subset g\i
M^s=M^{g\i sg}$, so that $M^{s'}=M^{g\i sg}$. Thus $s'=g\i sg$.
 $sg = gs'$.
\qed\end{demo}

Now the theorem follows from \nmb!{3.1}.
Indeed, by \nmb!{2.9} the group $G$ acts transitively
on the set of chambers. Assume that $gQ = Q$ for some $g \in G$. Then
$g \in \bigcap_{s\in S} P_{s} = \bigcap_{s\in S} P^+_{s_i}= 1$
by property \nmb!{3.1.1}.
\qed
\end{demo}

\begin{proclaim}{\nmb.{3.6}. Corollary}
A discrete group $G$ of isometries on a Riemannian manifold $M$ generated by
reflections is a quotient of a Coxeter group.
\end{proclaim}

\begin{demo}{Proof}
This follows by \nmb!{2.11}, or by \nmb!{2.13} and \nmb!{3.5}.
\qed\end{demo}

\subsection*{\nmb.{3.7}. Question }
Does there exist a discrete group of isometries which is generated by
reflections but is not a Coxeter group? If so, can one characterize
those which are Coxeter groups?

\begin{proclaim}{\nmb.{3.8}. Corollary}
Let $G$ be a reflection group on a complete connected
Riemannian manifold $M$ such that $G$ acts freely and
transitively on the set of all chambers, e.g., a Coxeter
manifold.  Let $C$ be a chamber.
Then we have:
\begin{enumerate}
\item[(\nmb:{1})]
  $C$ is the  Dirichlet domain associated with an interior point of $C$.
\item[(\nmb:{2})]
  Each chamber is convex  and its interior  consists of regular points.
\item[(\nmb:{3})]
  Any central hypersurface $H_{x_0,g.x_0}$ of a regular point $x_0$
  and $1\ne g\in G$ is a reflection hypersurface.
\item[(\nmb:{4})]
  $M_{\text{ reg}}=\bigcup_{g\in G}gC^o$.
\item[(\nmb:{5})]
  Let $F_1$ and $F_2$  be two walls of the chamber $C$ such that
  $F_1=g.F_2$ for some $g\in G$. Then $F_1=F_2$.
\item[(\nmb:{6})]
  The natural projection $\pi:M\to M/G$ induces a homeomorphism
  $C\to M/G$.
\end{enumerate}
\end{proclaim}

\begin{demo}{Proof}
{\nmb|{1}}
By lemma \nmb!{2.9} the chamber $C$ is a union of Dirichlet domains; but by
\nmb!{3.5} $G$ acts simply transitively on the set of chambers, thus $C$ is
just one Dirichlet domain, by \nmb!{2.3}.

{\nmb|{2}}
By {\nmb|{1}} and \nmb!{2.3.2} each chamber consists of regular points.
For convexity we have to show that any minimal geodesic arc between two
points in $C$ is contained in $C$. This follows from \cit!{2},~3.5.

{\nmb|{3}}
By {\nmb|{1}}, the union of all open chambers equals the union of all
open Dirichlet domains $D(x)$ for all regular points $x$. Thus also their
complements in $M$ are the same: The union of all reflection
hypersurfaces for $G$ in $M$ equals the union of all central
hypersurfaces with respect to some (each) regular point. Thus the
reflection hypersurfaces are exactly the central hypersurfaces
$H_{g.x_0,g'.x_0}$.

{\nmb|{4}}
If $x\in H_{x_0,g.x_0}$ then by {\nmb|{3}} the isotropy group of
$x$ is not trivial, so $x$ is not regular. Thus by
\nmb!{2.3.2} we have $M_{\text{reg}}=\bigcup_{g\in G}g.C^o$.

{\nmb|{5}} Let $F_1^o$ be the open interior of $F_1$ in
some central hypersurface $H$.
$F_1$ is contained in the intersection of exactly two chambers, namely
$F_1\subseteq C\cap h.C$, where $h$ is the reflection in the
hypersurface $H$. Also $F_2=g.F_1=g.C\cap g.h.C$, but one of
the two chambers must be $C$.
Thus $g=h\i$ is the reflection at $H$ and so $F_1=F_2$.

{\nmb|{6}} follows from {\nmb|{5}}
and from the fact that $G$ acts simply transitively on the set of all chambers.
\qed\end{demo}

\subsection*{\nmb.{3.9}}
Let $(M,G)$ be a connected Riemannian Coxeter manifold and let $C$
be a chamber. We denote by $W$ the set of walls of $C$ and by
$G(W)$ the free group, generated by involutive generators $r_F$
corresponding to all walls $F\in W$. Since $G$ is generated by
reflections with respect to walls in $W$, there is a natural
homomorphism $G(W) \to G$. We denote its kernel by $R$. We
define the normal subgroup $R_a$ of angular relations of
$G(W)$as follows:
\begin{quote}
  Let $F_i, F_j\in W$ be neighboring walls with non empty intersection $f$ 
  containing a codimension 2 submanifold, and
  let $F_i$ and $F_j$ have angle $\pi/n$ for a natural number $n$
  along some codimension 2 connected component of $f$, 
  then $(r_{F_i}r_{F_j})^{n}$ is a generator of $R_a$ in $G(W)$.
\end{quote}
We denote by $M_i$, $i=2,3$ the  complement in $M$ of the
union of codimension $\ge i$ intersections of reflection
hypersurfaces. Note that these intersections are totally 
geodesic  submanifolds  as  fixed  point  sets of finitely many
isometries. 

\begin{proclaim}{Theorem}
In this situation, the group $R_a$ of angular relations is a normal
subgroup of the group $R$ of all relations in $G$. 
Moreover, $\pi_1(M_3,x_0)=\pi_1(M,x_0)$ and $\pi_1(C^o,x_0)= \pi_1(C,x_0)$, and
we have the following exact sequences of groups:
\begin{gather*}
\{1\}\to \pi_1(C^o,x_0) *_e G(W) \to \pi_1(M_2,x_0) \to G(W) \to G \to \{1\} \\
\{1\}\to \pi_1(C,x_0) *_e G(W)/R_a\to \pi_1(M,x_0) \to G(W)/R_a \to G \to \{1\}
\end{gather*}
where for groups $H$ and $G$ the group $H*_e G$ is the kernel
of the projection $p_G:H*G\to G$ from the free product to $G$.
In particular, 
\begin{gather*}
\pi_1(M_2,x_0)/(\pi_1(C^o,x_0)*_e G(W))=R,\\
\pi_1(M,x_0)/(\pi_1(C,x_0)*_e G(W)/R_a)= R/R_a.
\end{gather*}
\end{proclaim}

\begin{demo}{Proof}
By \nmb!{3.8.6} the composition $C\to M\to M/G$ is a homeomorphism thus
$\pi_1(C,x_0)\to \pi_1(M,x_0)$ is injective.
By restriction $C\cap M_2 \to M_2 \to M_2/G$ is also a homeomorphism thus
$\pi_1(C\cap M_2,x_0)\to \pi_1(M_2,x_0)$ is injective.
By \nmb!{3.8.5} we have $\pi_1(C^o,x_0)= \pi_1(C,x_0)$ since a closed
curve in $C$ may be deformed into $C^o$.

Any element in $\pi_1(M,x_0)$ can be represented by a closed smooth curve
$c$ through $x_0$ in $M$ which we may assume to be transversal to all
intersections of walls. By dimension, $c$ lies in $M_2$ and first meets a
wall $F_1$ of $C$ transversally. Next it meets a wall $s_{F_1}(F_2)$ of
$s_{F_1}(C)$ transversally. And so on until it comes back to $x_0$. We assign
to $c$ the expression (word) $r_{F_1}r_{F_2}\dots r_{F_k}$ in $G(W)$. A homotopy
moving $c$ in $M_2$ just allows cancellations in this expression using
$r_F^2=1$. Replacing the $r_F$ in this expression by the corresponding
$s_F$ we get an element in the reflection group $G$ which maps $C$ to $C$
and thus is the identity, by theorem \nmb!{3.5}.

Let $f_i$ be a fixed curve from $x_0$ to $s_{F_i}(x_0)\in
s_{F_i}(C)$ hitting $F_i$ once transversally.
Any expression $r_{F_1}\dots r_{F_k}$ in $G(W)$ which maps to the
identity in $G$, is assigned to the closed curve in $M_2$ which 
first follows $f_1$ from
$x_0$ to $s_{F_1}(x_0)$, then 
$s_{F_1}\o f_2$ from
$s_{F_1}(x_0)$ to $s_{F_1}s_{F_2}(x_0)$, etc., until it ends
again in $x_0$.
Thus the sequence is exact at $G(W)$. 

A curve representing an element in $\pi_1(M_2,x_0)$ which is
transversal to walls 
can be described, up to `transversal'
homotopy, by a word
$c_0 r_{F_1} c_1 r_{F_2} c_2 \dots r_{F_k} c_k$ where: 
\begin{itemize}
  \item $c_i\in \pi_1(s_{F_1}s_{F_2}\dots
        s_{F_i}(C^o),s_{F_1}s_{F_2}\dots s_{F_i}(x_0)) \cong
        \pi_1(C^o,x_0)$, 
  \item $r_{F_i}$ stands for the curve 
        $s_{F_1}s_{F_2}\dots s_{F_{i-1}}(f_i)$. 
  \item $s_{F_1}s_{F_2}\dots s_{F_k}=e$ in $G$ since the curve
        is closed.
\end{itemize}
Thus the word describes a unique element of the free product
$\pi_1(C^o,x_0)* G(W)$ which is 
in the kernel of $\pi_1(C^o,x_0) * G(W)\to G$.
The curve in $\pi_1(M_2,x_0)$ maps to $e\in G(W)$ if and only
if the word above also satisfies
\begin{itemize}
  \item $r_{F_1}r_{F_2}\dots r_{F_k}=e$ in $G(W)$.
\end{itemize}
These are the elements of $\pi_1(C^o,x_0) *_e G(W)$.

So the first sequence is left exact, and surjectivity at $G$
follows from \nmb!{2.11}.

The second exact sequence follows from the first one:
any homotopy in $M$ between smooth curves in $M_2$ may be assumed to be
transversal to all intersections of reflection hypersurfaces of codimension
$\ge 2$. Then it avoids all intersection of codimension $\ge3$, so it lies
in $M_3$. Thus $\pi_1(M,x_0)=\pi_1(M_3,x_0)$. If the homotopy meets an
intersection $f=F_1\cap F_2$ transversely, moving the curve through $f$
means a cancellation in the expression assigned to the curve
which is given by the corresponding generator $(r_{F_1}r_{F_2})^n$ of $R_a$.
\qed\end{demo}

\begin{proclaim}{\nmb.{3.10}.  Theorem}
Let $(M,G)$ be a simply connected Riemannian Coxeter manifold and let $C$
be a chamber.
Then we have:
\begin{enumerate}
\item[(\nmb:{1})]        In terms of \nmb!{3.9} we have $R_a=R$. In other words,
  the relations \nmb!{2.11.1} and \nmb!{2.11.2} generate all relations of the
  Coxeter system $(G,S)$.
\item [(\nmb:{2})] The stabilizer $G_x$ of a point $x \in C$ is a
  finite Coxeter group generated by reflections with respect to the walls
  $F_{i_1}, \dots , F_{i_k}$ through $x$. Moreover, if $G_x$
  has no factor isomorphic to the dihedral group $D(m)$  for
  $m=5$ or $>6$, then  the angles between two walls through $x$ take
  values $\pi/n$ for $n = 2,3,4,6$.
\end{enumerate}
\end{proclaim}

For linear Coxeter groups this result was proved by Vinberg \cit!{32}.

\begin{demo}{Proof}
{\nmb|{1}}
This follows from $\pi_1(M,x_0)= R/R_a$ from \nmb!{3.9}.

{\nmb|{2}}
Let $g=s_{F_1}\dots s_{F_j}\in G_x$. Since any $h\in G$ preserves the union of
all reflection hypersurfaces, $g$ permutes the set of reflection hypersurfaces
through $x$. Thus $g(f)=f$ where $f$ is the connected component of
$F_{i_1}\cap \dots \cap F_{i_k}$ containing $x$. Then $C\cap gC\supseteq f$.

We shall use the method of proof of theorem \nmb!{3.9}.
Now choose a regular point $x_0\in C$ near $x$ and a curve $c_1$ in $M_2$
from $x_0$ to $gx_0$ which transverses the walls $F_j$, then
$s_{F_1}(F_{2})$, etc. Choose a second smooth curve $c_2$ in $M_2$ from $x_0$ to
$gx_0$ in $M_2$ which is near $x$ so that it intersects only walls through
$x$. Then we choose a homotopy in $M$ between $c_1$ and $c_2$ which we
may assume to be transversal to all codimension $\ge 2$ intersections of
reflection hypersurfaces. Then it is in $M_3$ and cuts intersections of two
reflection hypersurfaces transversely. Moving $c_1$ to $c_2$ via this
homotopy amounts to do angular cancellations (in $R_a$)
in the representation of $g$. Thus $g$ is represented also as a word in
reflections in hypersurfaces through $x$ according to the transversing of
$c_2$ of the corresponding walls.
\qed\end{demo}

\section*{\totoc\nmb0{4}. Riemannian manifolds with corners of Coxeter type }

\subsection*{\nmb.{4.1}. Manifolds with corners }
For more details see \cit!{25}, section 2.
A \idx{\it quadrant $Q\subset \mathbb R^n$ of index $k$}
is a subset of the form
$Q=\{x\in \mathbb R^n: l_1(x)\ge0,\dots,l_k(x)\ge0\}$ where
$l_1,\dots,l_k$ are independent linear functionals on $\mathbb R^n$.
If $x\in Q$ and exactly $j$ of the $l_i$ vanish on $x$ then
$x$ is called a \idx{\it corner of index $j$}.
For an open subset $U\subset Q$ a mapping
$f:U\to \mathbb R^p$ is called $C^r$ ($0\le r\le \infty$) if all partial
derivatives of $f$ of order $\le r$ exist and are continuous on $U$.
By the Whitney extension theorem this is the case if and only if $f$ can be
extended to a $C^r$ function $\tilde f:\tilde U\to \mathbb R^p$, where
$\tilde U\subset \mathbb R^n$ is open and $U=\tilde U\cap Q$.
If $f:U\to U'$ is a diffeomorphism between open subsets of quadrants
in $\mathbb R^n$ then the index of $x\in U$ equals the index of
$f(x)\in U'$.

A \idx{\it smooth manifold with corners} $M$ is defined in the usual way:
it is modelled on open subsets of quadrants in $\mathbb R^n$; a chart
on $M$ is a diffeomorphism $u:U\to u(U)$ from an open subset
$U\subset M$ onto an open subset $u(U)$ of a quadrant in $\mathbb R^n$,
where $n=\dim (M)$. The chart $(U,u,Q)$ is called centered on $x$ if
$u(x)=0$. A point $x\in M$ is called a \idx{\it corner of index $j$}
if there is a chart $(U,u,Q)$ of $M$ with $x\in U$ and and $u(x)$ a
corner of index $j$ in $Q$.

A subset $N\subset M$ is called a \idx{\it submanifold with corners}
of the manifold with corners $M$, if for any $y\in N$ there is a
chart $(U,u,Q)$ of $M$ centered at $y$ and there is a quadrant
$Q'\subseteq \mathbb R^k\subseteq \mathbb R^n$ such that $Q'\subseteq Q$
and $u(N\cap U)= u(U)\cap Q'$. A submanifold with corners $N$ of $M$
is called \idx{\it neat \ign{submanifold}} if the index in $N$ of
each $y\in N$ coincides with its index in $M$. Only neat submanfolds
have tubular neigborhoods.

Let us denote by $\partial^jM$ the set of all corners of index $j$ of
$M$. Note that $\partial^0M=M$.
Then each $\partial^jM$ is a submanifold without boundary of
$M$. Let $\partial M:=\bigcup_{j\ge 1}\partial^jM$.

Each closure (in $M$) of a connected component of
$\partial^jM$ is a submanifold with cornes of $M$ which is called a
{\it codimension $j$ face} of $M$; it is of dimension $n-j$.
A codimension 1 face is also called a {\it wall}.
A face is not neat. The set of all faces is a partially ordered set with
respect to inclusion.

The \idx{\it tangent bundle of a manifold with corners} $M$ is
constructed in the following way: Let $(U_\al,u_\al,Q_\al)$ be an
atlas of $M$. Then $TM$ is the quotient space of the disjoint union
$\bigsqcup_\al (\{\al\}\x U_\al\x \mathbb R^n)/\sim$
by the following equivalence relation: $(\al,x,v)\sim (\be,y,w)$ if
$x=y$ and $d(u_\be\o u_\al\i)(u_\al(x))v=w$. Then $\pi_M:TM\to M$ is
a smooth vector bundle, and the total space $TM$ is again a manifold
with corners: the corners are all in the base.

A tangent vector $X$ is called \idx{\it inner} (short for: not outer)
if there is a smooth curve $c:[0,1)\to M$ with $\dot c(0)= X$.
If $X\in T_xM$ and if $(U,u,Q)$ is a chart with $x\in U$, and if the
quadrant $Q$ is given by the independent linear functionals
$l_1,\dots,l_k$, and if $Tu(X)=(u(x),v)\in u(U)\x \mathbb R^n$,
then $X$ is inner if and only if the following holds: If
$l_i(u(x))=0$ then $l_i(v)\ge 0$, for all $i$.
Let us call the tangent vector \idx{\it strictly inner} if
$l_i(u(x))=0$ implies $l_i(v)>0$, for all $i$.
Let us denote the
space of all inner vectors by ${}^iTM\subset TM$. It is not a manifold
with corners any more. For example,
${}^iT[0,\infty) = \{(x,v):x\ge0, x=0\implies v\ge0\}$.

An inner vector field on $M$ is a smooth vector field $X:M\to TM$
whose values are all inner tangent vectors. By pasting local
solutions one can show that there exists a smooth open semiflow of $X$ in
the following sense: There is a set $W\subset \mathbb R\x M$ containing
$\{0\}\x M$ and $[0,\ep_x)\x \{x\}$ for some $\ep_x>0$ for each
$x\in M$ and a smooth mapping $\Fl^X:W\to M$ with $\Fl^X_0(x)=x$ and
$\frac{d}{dt}\Fl^X_t(x)=X(\Fl^X_t(x))$. But $\Fl^X_t$ is not even a
local diffeomorphism (it may map a corner to an interior point).

By a partition of unity argument one can show that there exists a
smooth vector field $Y$ on $M$ which is strictly inner, and one
may adapt it
in such a way that its flow $\Fl^Y_t$ is defined everywhere on $M$
for $0\le t\le\ep$ for $\ep>0$. Then $\Fl^Y_\ep$ maps $M$ into its
interior $M\setminus \partial M$. Thus: {\it Each manifold with
corners $M$ is a submanifold with corners of a manifold without
boundary of the same dimension.} See also \cit!{18}.

Let $X$ be a vector field on $M$ which is tangential to the boundary:
if $x\in\partial^jM$ then $X(x)\in T_x\partial^jM$ for all $j$. Then
there exists a local flow for $X$ for positive and for negative time;
the set $W\subset \mathbb R\x M$ is open.

\subsection*{\nmb.{4.2}. Equipment of a manifold with corners}
Let $M$ be an $n$-dimensional manifold with corners.
Consider a surjective mapping $s$ from the set $W$
of all walls (codimension 1 faces) of $M$ onto the set of generators $S$ of a
Coxeter system $(G,S)$ (see \nmb!{3.1}).
Any  face  $f$ of $M$ of codimension $k$ is the intersection of
$k$ many walls $W_1,\dots W_k$ (but not conversely). 
Then we extend the map $s$ to a map  $s$ from the set of faces of
$M$ into the set of Coxeter subgroups of $G$ as follows :
$$ s : f = F_1 \cap ... \cap F_k \mapsto s(f)= G(s(F_1), ... s(F_k))$$
where $G(s(F_1), ... s(F_k)$ is the subgroup of $G$
generated by  $s(F_1) , ... s(F_k)$.

The mapping $s$ is called a {\it Coxeter equipment} of $M$ by the Coxeter
system $(G,S)$, if $G(F)$ is a finite group for each
face of codimension $\ge 1$.
It follows that $s$ is an  partial order reversing
homomorphism of the poset of all
faces of $M$ into the poset of all Coxeter subgroups of the Coxeter system
$(G,S)$ if we also put $s(\emptyset)=G$. Note that $s(M)=\{1\}$.

\subsection*{\nmb.{4.3}. Riemannian manifolds with corners}
A \idx{\it Riemannian metric on a manifold with corners} $M$ is as
usual a smooth section $\ga:M\to S^2_+T^*M$. So it can be smoothly
extended to a Riemannian metric on a manifold without boundary of the
same dimension which contains $M$ as a submanifold with corners. If
the Riemannian metric has the property that each closure of a face is a
totally geodesic submanifold, then for each each inner tangent
vector $X_x\in {}^iT_xM$ the geodesic $t\mapsto\exp_x(tX_x)$ is defined
for small nonnegative $t$.

This can be expressed by the property of the
geodesic spray to be `inner' and `tangential' to all boundary strata
$\partial^jM$, see \cit!{25},~section~2. In detail: A vector
$\Xi\in TTM$ is called an \idx{\it inner tangent vector} to ${}^iTM$ if
there exists a smooth curve $c:[0,\ep)\to TM$ with
$\pi_{TM}(\Xi)=c(0)$, $c([0,\ep))\subset {}^iTM$, and $c'(0)=\Xi$. For
example, let $Q=\{x\in \mathbb R^n: l_1(x)\ge0,\dots,l_k(x)\ge0\}$ be a
quadrant and let $(x,u)\in {}^iTQ$. A vector $(x,u;v,w)\in T^2Q$ then
is inner to ${}^iTM$ if and only if:
\begin{enumerate}
\item[(\nmb:{1})] If $x$ is inner, so $u$ is arbitrary, then $(v,w)$ is
       arbitrary.
\item[(\nmb:{2})] If $l_i(x)=0$ and $l_i(u)>0$ then $l_i(v)\ge 0$ and $w$ is
       arbitrary.
\item[(\nmb:{3})] If $l_i(x)=0$ and $l_i(u)=0$ then $l_i(v)\ge 0$ and
       $l_i(w)\ge0$.
\end{enumerate}
Let us denote by ${}^iT^2M$ the set of all vectors which are inner to
${}^iTM$.
A \idx{\it spray} $S$ on the manifold with corners $M$ is a smooth
mapping $S:TM\to T^2M$ such that
\begin{enumerate}
\item[(\nmb:{4})] $T(\pi_M)\o S= \Id_{TM}$.
\item[(\nmb:{5})] $\pi_{TM}\o S = \Id_{TM}$.
\item[(\nmb:{6})] $T(m_t).S(X) = \tfrac 1t S(t.X)$ for $0\ne t\in \mathbb R$,
       where $m_t:TM\to TM$ is scalar multiplication by $t$.
\end{enumerate}
The spray is called \idx{\it inner\ign{ spray}} if
$S({}^iTM)\subset {}^iT^2M$ and it is called
\idx{\it tangential\ign{ spray}} if moreover $S$ is tantent to each
boundary stratum: $S(T\partial^jM)\subset T^2(\partial^jM)$.

If $\ga$ is a smooth Riemannian metric on the manifold with corners
$M$, then we may extend $\ga$ to a Riemannian metric $\tilde \ga$ on a
suitable open manifold $\tilde M$ of the same dimension which
contains $M$ as submanifold with boundary. We may compute the
geodesic (Levi-Civita) spray $\tilde S$ of $\tilde \ga$ and restrict it
again to $TM$. This spray is an inner tangential spray if and only
if in $(M,\ga)$ all closures of faces are totally geodesic submanifolds,
and we have $\exp = \pi_M\o \Fl^S_1$.

Thus we conclude (see also \cit!{25},~2.10):

\begin{proclaim}{Lemma} {\rm \cit!{25},~2.10}
Let $\ga$ be a Riemannian metric on a manifold with corners $M$
such that all faces are totally geodesic.
Then there exists a suitable open neighborhood $V$ of the zero
section in $T^iM$ such that the geodesic exponential mapping
$\exp:V\to M$ is defined. If $V$ is small enough then $\exp$ has the
following properties:
\begin{enumerate}
\item $\exp(0_x)=x$ for all $x\in M$.
\item $\exp_x:V_x:= V\cap T^i_xM \to M$ is a diffeomorphism of $V_x$
      onto an open neighborhood $W_x$ of $x$ in $M$.
\item $V_x$ is the intersection of an open ball
       $B_x\subset (T_xM,\ga_x)$ with a quadrant $Q_x\subset T_xM$.
\item The mapping $(\pi_M,\exp):V\to M\x M$ is a diffeomorphism onto an
      open neighborhood of the diagonal in $M\x M$.
\item $\exp$ restricts to the exponential mapping of the induced
      Riemannian metric on each closure of a face.
\end{enumerate}
\end{proclaim}

\subsection*{\nmb.{4.4}. Riemannian chambers and their Coxeter equipment }
An \idx{\it Riemannian chamber} is a manifold with corners
$C$ with a Riemannian metric $\ga$ such that each face is totally
geodesic and such that the following two conditions {\nmb|{1}}
and {\nmb|{2}} are satisfied.
\begin{enumerate}
\item[(\nmb:{1})] The angle between neighboring walls $W_i$ and $W_j$
is  a  constant of the form $\pi/n_{ij}$ for $n_{ij}\in \mathbb
N$  along  any  codimension  2  connected component of $W_i\cap
W_j$. 
\end{enumerate}
Let $V\subset T^iC$ be small as in \nmb!{4.3}. Then
$\exp_x:V_x=V\cap T^i_xC\to W_x\subset C$ is a
diffeomorphism. Recall from \nmb!{4.3} that $V_x$  is the
intersection of an open ball $B_x$ in $(T_xC,\ga_x)$ with a
quadrant $Q_x$ whose walls contain the inverse images under $\exp_x$
of the closed walls of $C$ containing $x$. The angles between the
hyperplanes $T_xW_i$ and $T_xW_j$ in the Euclidean space $(T_xC,\ga_x)$ are exactly
$\pi/n_{ij}$, by {\nmb|{1}}.
By \cit!{32}, theorem 1, this equivalent to the fact that
the group $G_x\subset O(T_xC,g_x)$
generated by the reflections in the hyperplanes $T_xW_i$ is a finite
Coxeter group with fundamental Weyl chamber $\mathbb R_{\ge0}.V_x$.

Consider the pullback Riemannian
metric $(\exp_x|V_x)^*\ga$ on $V_x$. Now we can formulate the second
condition:
\begin{enumerate}
\item[(\nmb:{2})] If we extend the Riemannian metric $(\exp_x|V_x)^*\ga$ on
       $V_x$ to the ball $B_x=G_x.V_x$ by using the elements of $G_x$ as
       isometries, then the resulting $G_x$-invariant Riemannian
       metric $\tilde \ga_{B_x}$ on $B_x$ is smooth.
\end{enumerate}
If $G$ is a discrete group of isometries of a complete Riemannian
manifold $(M,\ga)$ which is generated by disecting reflections, and if
a chamber $C$ is also a Dirichlet domain, then
obviously $(C, \ga)$ is a Riemannian chamber.

\begin{proclaim}{Proposition}
Any Riemannian chamber $C$ carries a universal Coxeter
equipment.
\end{proclaim}

\begin{demo}{Proof}
Let $\{W_i\}$ be the set of all walls of $C$.
For each wall $W_i$ of
$C$ we take a generator $s_i$.
Then let $G$ be the group generated by all $s_i$, with relations
$(s_is_j)^{n_{ij}}=1$, whenever $W_i\cap W_j\ne \emptyset$ and where the angle
between $W_i$ and $W_j$ is $\pi/n_{ij}$.
Then $G$ is a Coxeter group with Coxeter system $(G,\{s_i\})$.
For each $x\in C$ we
constructed in \nmb!{4.4} a linear Coxeter group
$G_x\subset O(T_xC,g_x)$ which is generated by those $s_i$ for
which $x\in W_i$. Obviously, $G_x$ is a finite subgroup of $G$.
Moreover, let $F=W_{1}\cap \dots W_k$ be a nonempty face. Then $G(F)$
is generated by the reflections $s_1,\dots,s_k$ which satisfy pairwise
$(s_is_j)^{n_{ij}}=1$ for $2\le n_{ij}<\infty$. Thus $G(F)$ is finite
for each nonempty face.
\qed\end{demo}

The Coxeter equipment constructed in this proposition is called universal
since the mapping $s$ is injective. Other Coxeter equipments are possible,
if different walls are mapped to the same generator in such a way, that the
isotropy group of each face $F$ stay isomorphic to $G(F)$ as above,
and the full group is still a Coxeter group.

Thus we say that a Coxeter equipment $s$ of the Riemannian chamber $C$ is
admissible, if for any two different walls $W_i$ and $W_j$ with nonempty intersection
the element $s(W_i)s(W_j)$ has order exactly $n_{ij}$ in $G$, where the
angle between $W_i$ and $W_j$ is $\pi/n_{ij}$. The right hand side of
fig.\ ~\fig!{2} gives an example of a not universal equipment.

\subsection*{\nmb.{4.5}. The Coxeter Riemannian manifolds associated with a
Riemannian chamber }
Note that by \nmb!{2.11} and \nmb!{3.9} the Weyl chamber $C$ of
a Coxeter $G$-manifold $M$ has the natural
structure of a Riemannian chamber with the admissible
equipment $s : F \mapsto s(F) = \langle s \in S: M^s \supset
F\rangle$. In the non-disecting case this is not true: In \nmb!{2.10} 
the chamber of the non-disecting reflection on $\mathbb R\mathbb P^2$ equals
$\mathbb R\mathbb P^2$ and the generating reflection is not associated to a wall since 
$\mathbb R\mathbb P^2$ has no boundary. 

For Coxeter manifolds the converse statement is also true as the
following theorem shows.

\begin{proclaim}{Theorem}
Let $C$ be a Riemannian chamber.

Then to each admissible Coxeter equipment $G$ of $C$
there exists a
smooth Riemannian manifold $\mathcal U(G,C)$ without boundary and a discrete
subgroup $G$ of isometries which is generated by 
reflections such that $C$ is isometric to a chamber of $M$ which
is also a Dirichlet domain.

If $C$ is connected then also $M$ is connected.
If the equipment $G=G_{\text{univ}}$ is the universal one
then $G$ is generated by disecting reflections and 
$\pi_1(\mathcal U(G,C))=\pi_1(C)*_e G_{\text{univ}}$. In general
we have an exact sequence:
\begin{displaymath}
\{1\}\to \pi_1(C)*_e G_{\text{univ}} \to \pi_1(\mathcal U(G,C)) 
\to G_{\text{univ}} \to G \to \{1\}.
\end{displaymath}
\end{proclaim}

\begin{demo}{Proof} We use first the universal equipment.
Let $\{F_i\}$ be the set of all closures of walls of $C$. We
construct first the group $G$, as follows. For each wall $F_i$ of
$C$ we take a generator $s_i$ of $G$.
Then $G=G_{\text{univ}}$ is the group generated by all $s_i$ and with relations
$(s_is_j)^{n_{ij}}=1$, when $F_i\cap F_j\ne \emptyset$ and where the angle
between $F_i$ and $F_j$ is $\pi/n_{ij}$.
For each $x\in C$ we
constructed in \nmb!{4.4} a linear Coxeter group
$G_x\subset O(T_xC,g_x)$ which is generated by those $s_i$ for
which $x\in F_i$. Obviously, $G_x$ is a subgroup of $G$.

Now we construct $M=\mathcal U(G,C)$ as topological space by putting
$\mathcal U(G,C):= G\x C/\sim$ where
\begin{gather*}
(g.s_i,x)\sim (g, s_i(x))=(g, x) \text{  for }x\in F_i,\quad
\text{or equivalently}\\
(g,x)\sim (h,y)\Longleftrightarrow x=y\text{ and } g\i h\in G_x.
\end{gather*}
So $\mathcal U(G,C)$ is a quotient of the disjoint union of $|G|$
copies of $C$ which are glued together only along walls.

We construct an atlas for $\mathcal U(G,C)$ as follows, using the arguments
from \nmb!{4.4}. For a corner $x$ of $C$ consider the Riemannian
metric $\tilde \ga_{B_x}$ on the open ball $B_x\subset T_xC$ which is
smooth by condition \nmb!{4.4.2},
and the smooth exponential mapping
$\exp_x:V_x= B_x\cap Q_x \to W_x\subset C$.
We extend it to a $G_x$-equivariant homeomorphism
$\tilde{\exp_x}$ from $B_x$ to the open neighborhood
$U_x=\bigcup_{g\in G_x}(\{g\}\x W_x)$ of $x$ in $M$ by putting
$\tilde{\exp_x}(g.X)=(g,\exp_x(X))$ for $X\in V_x$ and $g\in G_x$.
Then $(U_x,u_x:=\tilde{\exp_x}\i:U_x\to B_x\in T_xC)$ is a chart
on $M$.

If $x\in C$ is a regular point we use the inverse of the
exponential mapping on such a small neigborhood of 0 in $T_xC$
that its image does not meet any wall. These charts we the distribute
from $C=\{\Id\}\x C$ to the whole of $M$ by using the
transformations from $g$.

We claim that this gives a smooth atlas for $\mathcal U(G,C)$: Suppose that $x$ and
$y$ are corners of $C$ such that $W_x\cap W_y\ne \emptyset$. We
have to show that $u_x\o u_y\i$ is smooth. We may assume that
$y\in W_x$ since we may connect $x$ and $y$ by finitely many chart
changings with this property. But then this is a chart change of
exponential mappings at different base points of the smooth
Riemannian metric in $B_x\subset T_xC$.

Finally, $G$ acts on the smooth manifold $\mathcal U(G,C)$ by construction:
$g.(g_1,x)=(gg_1,x)$, and it consists of isometries. 
By construction $G$ acts freely and transitively on the set of
all chambers of $\mathcal U(C,G)$.
We claim
that the generators $s_i$ of $G$ are disecting. Suppose for
contradiction that a generator $s$ is not disecting. Choose regular point
$x_0\in C^o$ and a smooth curve $c$ in $\mathcal U(C,G)\setminus
\mathcal U(C,G)^s$ from $x_0$ to $s.x_0$ which is transversal
to all intersections of reflection hypersurfaces. Then $c$
passes from $C$ to a neighbor $s_{i_1}C$, then to a neighbor
$s_{i_1}s_{i_2}C$ of $s_{i_1}C$, and so on, till it reaches the
chamber $s_{i_1}\dots s_{i_k}C=sC$ containing $s.x_0$. None of
the $s_{i_j}$ equals $s$ since $c$ does not meet $\mathcal
U(C,G)^s$.
Since $G$ acts freely and transitively on the set
of chambers we have $s=s_{i_1}\dots s_{i_k}$ in $G$, a
contradiction. (WHY??)

Finally, for a general admissible equipment we have a normal subgroup
$R\subset G$ of further relations which by the description of an admissible
equipment acts freely and discretely on the universal $M$ which thus is a
covering of the resulting manifold. The statement on fundamental groups
follows from \nmb!{3.9}.
\qed\end{demo}

\subsection*{\nmb.{4.6}. Remark }
We can also consider manifold with corners $C$ with a smooth
Riemannian metric $g$ which satisfies only condition
\nmb!{4.4.1}. Then we can construct a topological
manifold $M$ which is smooth off the union of all reflection
hypersurfaces, with a Riemannian metric which is only continuous along
the the reflection hypersurfaces, in general. It might be worthwile to
study this object.

\begin{proclaim}{\nmb.{4.7}. Theorem}
Let $C$ be a manifold with corners with a Coxeter equipment
$s:W\mapsto s(W)\in S$ where $(G,S)$ is a Coxeter system.

Then there exists a Riemannian metric $\ga$ such that $(C,\ga)$ is a
Riemannian chamber and $s$ is an admissible equipment for it.
\end{proclaim}

\begin{demo}{Proof} We construct the metric inductively starting from faces
which are manifolds without boundary.
On each such face $F$ we put an
arbitrary Riemannian metric $\ga_F$.

Now let $F$ be a face which contains corners of index (in $F$) at most 1,
i.e., $F$ is a
manifold with boundary $\p F$ which is a disjoint union of faces $F_1$ without
boundary.
Along each boundary component $F_1$ of $F$ we consider an open collar
$F_1\x [0,1)\subset F$
and extend the metric by $\ga(x,t)= dt^2+\ga_{F_1}(x)$ where $x\in F_1$ and $t$
is the coordinate function on $[0,1)$.
With a partition of unity we may extend this metric to
the whole of $F$ in such a way that near each $F_1$ it is not changed.
Note that $F_1$ is totally geodesic in $F$, and that
the metric is constant in the direction $t$ normal to $F_1$.

Now let $F$ be a face which contains corners of index (in $F$) at most 2,
i.e., $\p F$ contains walls $F_1^i$ of $F$ which are manifolds with
boundary. We already defined Riemannian metrics $\ga_{F_1^i}$
on $F_1^i$.
If $F_2$ is a boundary component of $F_1^1\cap F_1^2$ we consider an open
tubular wedge neigbourhood $F_2\x D$ of $F_2$ in $F$
with the following property. Each fiber $\{x\}\x D$ intersects $F_1^i$
exactly in the fiber $\{x\}\x [0,1)$ of that collar of $F_2$ in $F_1^i$ for
each $x\in F_2$ which was used above to construct
the Riemannian metric on
$F_1^i$. The fiber $D\subset \mathbb R^2$ is an open 0-neighborhood in a quadrant
with angle $\al_F(F_1^1,F_1^2)$ as in fig.\  \fig!{9}.
\begin{figure}[h]
\psfrag{F11}{$F^1_1$}
\psfrag{F12}{$F^2_1$}
\psfrag{F2}{$F_2$}
\psfrag{D}{$D$}
\psfrag{a}{$\al$}
\centerline{\includegraphics[width=8cm]{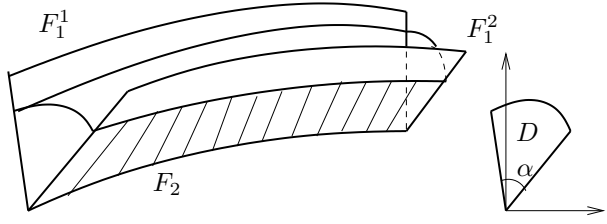}}
\caption{{\fig.{9}}
The tubular wedge neigborhood $F_2\x D$ and its fiber $D$.
}
\end{figure}\noindent
Here $\al_F(F_1^1,F_1^2)$ is determined by the Coxeter equipment:
If in terms of walls $W_i$ of $C$ we have
\begin{gather*}
F=W_{i_1}\cap\dots\cap W_{i_{n-2}},\\
F_1^1=W_{i_1}\cap\dots\cap W_{i_{n-2}}\cap W_{i_{n-1}},\qquad
F_1^2=W_{i_1}\cap\dots\cap W_{i_{n-2}}\cap W_{i_{n}},\\
F_2=W_{i_1}\cap\dots\cap W_{i_{n-2}}\cap W_{i_{n-1}}\cap W_{i_n},
\end{gather*}
then $\al_F(F_1^1,F_1^2)$ is determined by the (finite) Coxeter
system $(G(F_2)=G(S'),S')$
where $S'=\{s_{i_1},\dots,s_{i_n}\}$, by considering
the angle between and in the corresponding faces in the
Weyl chamber of $(G(S'),S')$,
as decribed in \nmb!{3.3}.

We now put the product metric $\ga_D(u)+\ga_{F_2}(x)$ for $(x,u)\in F_2\x D$
on the tubular wedge neigbourhood $F_2\x D$, where $\ga_D$ is the standard
Euclidean metric on $\mathbb R^2$ restricted to $D$. This gives a metric on
$F_2\x D$ which induces the already constructed metric $\ga_{F_1^i}$ on the
intersection with $F_1^i$ since $F_2\x D$ intersects $F_1^i$ in the collar
used to construct $\ga_{F_1^i}$. Moreover $F_2$ and the parts of $F_1^i$
are totally geodesic, and the metric is constant in directions normal to
any relevant face, near that face.

We do this construction near any face of codimension 2 of $F$.
Then we use a collar $(F_1^i\setminus \p F_1^i)\x [0,1)$ of the interior of
the face $F_1^i$ in $F$ such that the fiber near any $F_2$ coincides with the
normal geodesic in $F_2\x D$ in the metric constructed there.
Put the metric $dt^2+\ga_{F_1^i}(x)$ for
$(x,t)\in (F_1^i\setminus \p F_1^i)\x [0,1)$ on this collar, and use a
partition of unity on the union of all these collars and the wedge
neighborhoods which is constant in the normals near any face to glue the
metrics in such a way that the resulting metric is constant in the normal
directions near any face and each face is totally geodesic. With another
partition of unity we extend this metric into the interior of $F$ and not
changing it near any face.

We proceed inductively. We assume that we have already constructed in this
way metrics on each face which consists of corners of index $\ge k$ in $C$
and consider now a face $F$ which consists of corners of index $\ge k-1$ in
$C$. Then the boundary $\p F$ is a union of faces where we alredy
constructed the metric. Let $F_k$ be a minimal face in $\p F$, i.e., $F_k$
does not contain any other face. Then $F_k$ is a manifold without boundary
where we already have a metric $\ga_{F_k}$. Moreover $F_k$ is the
transversal intersection of $k$ walls $F_1^1,\dots,F_1^k$ of $F$, where $k$
is the codimension of $F_k$ in $F$. We then choose a tubular wedge
neighbourhood $F_k\x D^k$ of $F_k$ in $F$ which intersects fiber
respectingly each intersection of $k-1$ of the walls $F_1^1,\dots,F_1^k$ of
$F$ in the tubular wedge neighborhood which was used previously to
construct the metric $\ga_{F_1^i}$ on each of the walls. Here $D^k$ is an
open 0-neighborhood in a quadrant in $\mathbb R^k$ with walls whose angles
$\al_F(F_1^i, F_1^j)$ are determined by the
Coxeter equipment as described above.
We now put the metric $\ga_{D^k}(u)+\ga_{F_k}(x)$ for $(x,u)\in F_k\x D^k$
on the tubular wedge neigbourhood $F_k\x D^k$, where $\ga_{D^k}$ is the standard
Euclidean metric on $\mathbb R^k$ restricted to $D^k$. This gives a metric on
$F_k\x D^k$ which induces the already constructed metric $\ga_{F_1^i}$ on the
intersection with $F_1^i$ since $F_k\x D^k$ intersects $F_1^i$ in the
tubular wedge neighborhood
used to construct $\ga_{F_1^i}$. Moreover $F_k$ and the parts of $F_1^i$
are totally geodesic, and the metric is constant in directions normal to
any face near that face.

We do this construction near any minimal face of $F$.
Then we use a collar $(F_1^i\setminus \p F_1^i)\x [0,1)$ of the interior of
the face $F_1^i$ in $F$ such that the fiber near any minimal face $F_l$
coincides with the normal geodesic in $F_l\x D^l$ in the metric constructed there.
Put the metric $dt^2+\ga_{F_1^i}(x)$ for
$(x,t)\in (F_1^i\setminus \p F_1^i)\x [0,1)$ on this collar, and use a
partition of unity on the union of all these collars and the wedge
neighborhoods which is constant in the normals near any face, to glue the
metrics in such a way that the resulting metric is constant in the normal
directions near any face and each face is totally geodesic. With another
partition of unity we extend this metric into the interior of $F$.

Eventually we exhaust each connected component of $C$.
\qed\end{demo}

\begin{proclaim}{\nmb.{4.8}. Proposition}
Let $C$ be a manifold with corners with a Coxeter equipment
$s:W\mapsto s(W)\in S$ where $(G,S)$ is a Coxeter system.
Let $\ga$ and $\ga'$ be 
two Riemannian metrics on $C$ such that $(C,\ga)$ and $(C,\ga')$ are both
Riemannian chambers and $s$ is an admissible equipment for both.

Then the smooth manifolds $\mathcal U(G,C,\ga)$ and $\mathcal U(G,C,\ga')$
constructed via \nmb!{4.5} are diffeomorphic.
\end{proclaim}

\begin{demo}{Proof}
Since the construction as a topological space described in the proof of
\nmb!{4.5} depends only on the equipment, the two manifolds are
canonically homeomorphic. 
For a corner $x\in C$ let $u_x:U_x\to B_x\subset T_xC$ and 
$u_x':U_x'\to B_x\subset T_xC$ be two charts as described in the proof of 
\nmb!{4.5} for the two Riemannian metrics $\ga$ and $\ga'$. But then the
chart change $u_x'\o u_x\i$, considered in a manifold without boundary
which contains $C$ as a submanifold with corners (see \nmb!{4.1}), consists
of the exponential mapping of the extended Riemannian metric $\tilde\ga$ followed by the
inverse of the exponential mapping of $\tilde\ga'$, which is obviously
smooth. Thus the canonical homeomorphism between $\mathcal U(G,C,\ga)$ and 
$\mathcal U(G,C,\ga')$ is a diffeomorphism.
\qed\end{demo}

\section*{\totoc\nmb0{5}. Orbifolds
}

\subsection*{\nmb.{5.1}. Smooth orbifolds  }
We recall the definition of orbifold.
Let $X$ be  a second countable Hausdorff
space. An atlas of a smooth $n$-dimensional orbifold (or $V$-manifold) on $X$
is a family $\{U_i\}_{i\in I}$ of open sets that satisfy:
\begin{enumerate}
\item $\{U_i\}_{i\in I}$ is an open cover of X.
\item For each $i\in I$ a local uniformizing system consisting of
     a triple $\{\tilde U_i,G_i,\ph_i\}$, where $\tilde U_i$ is a
     connected open subset of $\mathbb R^n$ containing the origin,
     $G_i$ is a finite
     group of diffeomorphisms acting effectively and properly on $\tilde U_i$, and
     $\ph_i:\tilde U_i\to U_i$ is a continuous map of $\tilde U_i$ onto $U_i$
     such that $\ph_i\circ g=\ph_i$ for all $g\in G_i$ and
     the induced map of $\tilde U_i/G_i$ onto $U_i$ is a homeomorphism.
     The finite group $G_i$ is called a local uniformizing group.
     \item Given $\tilde x_i\in\tilde U_i$ and $\tilde x_j\in\tilde U_j$ such that
     $\ph_i(\tilde x_i)=\ph_j(\tilde x_j)$, there is a diffeomorphism
     $\ph_{ij}:\tilde V_j\to\tilde V_i$ from a neighborhood
     $\tilde V_i\subset\tilde U_i$ of $\tilde x_i$ onto a neighborhood
     $\tilde V_j\subset\tilde U_j$ of $\tilde x_j$ such that
     $\ph_i=\ph_j\circ \ph_{ji}$.
\end{enumerate}
Two atlases are equivalent if their union is again an atlas of a smooth
orbifold on $X$. An orbifold is the space $X$ with an equivalence class
of atlaces of smooth orbifolds on $X$.

\begin{proclaim}{Proposition} \cit!{31}
If $M$ is an $n$-dimensional smooth manifold and $G$ is a group acting
smoothly and discretely 
on $M$, then $X=M/G$ has a structure of orbifold.
\end{proclaim}

\begin{demo}{Proof} Let $x\in X$. Choose $\tilde x\in M$ projecting to
$x$, and denote by $G_x$ the isotropy group of $\tilde x$.
Choose a neigborhood of $\tilde U_x$ invariant by $G_x$ and disjoint
from $g(U_x)$ for all $g\in G\setminus G_x$ such that there
is a local chart $k:U_x\to\tilde U_x\subset\mathbb R^n$ on $M$ with $k(x)=0$.
We take
$(\tilde U_x,G_x,\ph_x)$, where $\ph_x$ is a composition of
$k^{-1}$ with the projection $\tilde U\to\tilde U/G_x$,
for a local uniformizing system. It is easily checked such local uniformizing
systems form an atlas of a smooth $n$-dimensional orbifold on $M/G$.
\qed\end{demo}

In the definition of atlas of a smooth orbifold on $X$ we can always
take the finite subgroups $G_i$ to be subgroups of the
orthogonal group $O(n)$ acting naturally on $\mathbb R^n$.
Condition (3) implies that for each
$g_i\in G_i$ there exists $g_j\in G_j$ such that
$\ph_{ji}\circ g_i=g_j\circ \ph_{ji}$.

Let $\{\tilde U_i,G_i,\ph_i\}$ be a unifomizing system such that
$\tilde U_i$ contains the origin, the group $G_i$ is a subgroup of
$O(n)$, and $x=\ph_i(0)$. Then the group $G_x=G_i$ is
independent of
the uniformizing system $\{\tilde U_i,G_i,\ph_i\}$. More precisely,
this group is defined up to
isomorphism and its action on $\mathbb R^n$ is defined up to isomorphism
as well. The point $x\in X$ is called regular if the corresponding group
$G_x$ is trivial and otherwise singular.

\subsection*{\nmb.{5.2}. Reconstruction of the orbifold structure from
the structure sheaf }
Let again $\{\tilde U_i,G_i,\varphi_i\}$ be a unifomizing system
such that $\tilde U_i$ contains the origin, the group $G_i$ is
a subgroup of $O(n)$, and $x=\varphi_i(0)$. Then there is a representation
$\rho:G_i\to O(n)$, a ball $B$ in $\mathbb R^n$ centered at the
origin, and a map $\varphi:B\to X$ such that $\varphi(0)=x$ and
$\{B,G_i,\varphi\}$ is a uniformizing system of the orbifold
$X$.

A function $f:U_i\to\mathbb R$ is called smooth
if $f\circ\varphi_i$ is a smooth function on $\tilde U$. The germs
of smooth functions on $X$ define a sheaf $S_X$ on $X$.

\begin{proclaim}{\nmb.{5.3}. Definition} Let $X$ and $\tilde X$ be two
smooth orbifolds.
The orbifold $\tilde X$ is called a covering orbifold for $X$
with a projection $p:\tilde X\to X$ if $p$ is a continuous map of underlying
topological spaces and  each point $x\in X$ has a neighborhood
$U=\tilde U/G$ (where $\tilde U$ is an open subset of $\mathbb R^n$)
for which each component $V_i$ of $p^{-1}(U)$ is isomorphic to
$\tilde U/G_i$, where $G_i\subset G$ is some subgroup. The
above isomorphisms $U=\tilde U/G$ and $V_i=\tilde U/G_i$ must
respect the projections.
\end{proclaim}

Note that the projection $p$ in the above definition is not a cover of
underlying topological spaces.

Hereafter we suppose that all orbifolds and their covering orbifolds
are connected.

\begin{proclaim}{\nmb.{5.4}. Theorem} \cit!{31}
An orbifold $X$ has a universal covering
orbifold $p:\tilde X\to X$. More precisely, if $x\in X$,
$\tilde x\in \tilde X$ are regular points and $p(\tilde x)=x$,
for any other covering orbifold $p':\tilde X'\to X$ and
$\tilde x'\in\tilde X'$ such that $p'(\tilde x')=x$ there is a cover
$q:\tilde X\to\tilde X'$ such that $p=p'\circ q$ and $q(\tilde x)=\tilde x'$.
For any points $\tilde x,\tilde x'\in p^{-1}(x)$ there is a deck
transformation of $\tilde X$ taking $\tilde x$ to $\tilde x'$.
\end{proclaim}

Suppose $\rho:G\to O(n)$ is a representation
of a finite group $G$, $\mathbb R^n/G$ is the
corresponding orbifold, and $S_{\mathbb R^n/G}$ is the corresponding
sheaf. By the Hilbert theorem the ring $\mathbb R[\mathbb R^n]^{G}$ is
finitely generated. Let $\sigma^1,\dots,\sigma^m$ be a system of homogeneous
generators of $\mathbb R[\mathbb R^n]^{G}$ and $y^1,\dots,y^m$ the
corresponding functions on $\mathbb R^n/G$. Consider the map
$\sigma=(\sigma^1,\dots,\sigma^m):\mathbb R^n\to\mathbb R^m$ called the orbit
map. It is known [2] that the map $\sigma$ induces a homeomorphism between
$\sigma(\mathbb R^n)$ and the orbit space $\mathbb R^n/G$ which establishes
an isomorphism between the restriction of the sheaf $C^{\infty}_m$
of smooth functions on $\mathbb R^m$ to $\sigma(\mathbb R^n)$ and the sheaf
$S_{\mathbb R^n/G}$.

It is clear that for each orbifold $X$ and $x\in X$ there is
a neghborhood $U_x$ and a representation
$\rho:G_x\to O(n)$ such that the restriction of $S_X$ to $U_x$ is
isomorphic to the restriction of the sheaf $S_{\mathbb R^n/G_x}$ to some
ball centered at the origin.

For a representation $\rho:G\to O(n)$ a diffeomorphism of the
orbit space $\mathbb R^n/G$ is an automorphism of the sheaf
$S_{\mathbb R^n/G}$ by definition. Let
$f:\mathbb R^n/G\to\mathbb R^n/G$ be a
diffeomorphism and $h^1,\dots,h^m$ a system of generators of
$S_{\mathbb R^n/G}$. Then $f$ is uniquely defined by
the images of generators
$h^i$ and these images are the generators of $S_{\mathbb R^n/G}$
again.
Denote by $R$ the set of all reflections contained in $G$ and by
$A(G,R)$ the set of all automorphisms of the group $G$ which
preserves the set $R$.

\begin{proclaim}{\nmb.{5.5}. Theorem}
\cit!{23} For each diffeomorphism $f$ of the orbit space
$\mathbb R^n/G$ there is a smooth lift $F:\mathbb R^n\to\mathbb R^n$.
For each such lift $F$ there is an automorphism $a\in A(G,R)$
such that for all $g\in G$ and $x\in \mathbb R^n/G$ we have 
$F(g x)=a(g)F(x)$.

The local version of this theorem is also true, i.e. if $B$ is a
ball in $\mathbb R^n$ centered at the origin and $f$ is a
diffeomorphism of the sheaf $S_{B/G}$, then there is a smooth lift
$F:B\to B$ with the same property as above.
\end{proclaim}

\begin{proclaim}{\nmb.{5.6}. Theorem}
An orbifold $X$ is defined
uniquely by its sheaf $S_X$.
\end{proclaim}

\begin{demo}{Proof} Note that for a regular point $x\in X$ the ring $S_X(x)$ of the
germs of $S_X$ at $x$ is isomorphic to the ring of germs at 0 of smooth
functions on $\mathbb R^n$.
Then the dimension of the orbifold $X$ is defined by the sheaf $S_X$.
Next note that if $\rho:G\to O(n)$ is a representation of a finite
group $G$, then the group preserving all smooth $G$-invariant
functions on $\mathbb R^n$ coincides with $\rho(G)$.
If this group is infinite there is a regular point with
non trivial stabilizer, which is impossible.  The result then follows from
the fact that the order of $G$ equals the cardinality
of a regular orbit.

It is sufficient to prove that for each a finite group $G$, a
representation
$G\to O(n)$, a ball $B$ in $\mathbb R^n$, and the
map $\varphi:B\to X$ which induces an isomorphism of the
sheaf $S_{B/G}$ and the restriction $S_U$ of the sheaf $S_X$ to some
open subset $U$ of $X$,
$\{B,G,\varphi\}$ is a uniformizing system on $X$.

Let $\{B_1,G_1,\varphi_1\}$ be such a uniformizing system,
corresponding to the representation $\rho_1:G_1\to O(n)$,
$\varphi(0)=x$,
and $\{B_2,G_2,\varphi_2\}$  a uniformizing system of the
orbifold $X$ which is induced by some representation
$\rho_2:G_2\to O(n)$  such that $\varphi_2(0)=x$. We may assume that
$B_1=B_2=B$ and
$\varphi_1(B)=\varphi_2(B)=U$.
Then the rings of functions on $B$ which are compositions of
$\varphi_1$ and $\varphi_2$ with the sections of $S_X$ on $U$ coincides.
By the above remark $\rho_1(G_1)=\rho_2(G_2)=G$.

For $i=1,2$ denote by $\bar\ph_i$ the diffeomorphism $B/G\to U$ induced by
$\ph_i$. Then $\bar\ph_2\i\o\bar\ph_1$ is a diffeomorphism of $B/G$.
By Theorem \nmb!{5.5} there is a
smooth lift $B\to B$ of this diffeomorphism. But this means
that $\{B_1,G_1,\varphi_1\}$ is a uniformizing system of
the orbifold $X$.
\qed\end{demo}

\begin{proclaim}{\nmb.{5.7}. Corollary}
Let a group $G$ acts discretely on a smooth simply
connected manifold $M$ and $S_X$ the corresponding sheaf on $X=M/G$.
Then $M$ is a universal covering orbifold for $X$.
\end{proclaim}

\begin{demo}{Proof} Evidently manifold $M$ is a covering orbifold for $X$.
If $\tilde X$ is universal covering orbifold for $X$, then there is
a cover $q:\tilde X\to M$. By the definition of cover $\tilde X$ should
be a manifold and $q$ a cover of manifilds. Therefore $q$ is a diffeomorphism.
\qed\end{demo}
Theorems \nmb!{5.6} and \nmb!{5.7} imply the following statement.

\begin{proclaim}{\nmb.{5.8}. Corollary}
Let a group $G$ act discretely on a smooth simply
connected manifold $M$ and $S_X$ the corresponding sheaf on $X=M/G$. Then
each diffeomorphism of the orbit space $X$, i.e. an automorphism of the
sheaf $S_X$ has a smooth lift to $M$.
\end{proclaim}

\section*{\totoc\nmb0{6}. Coxeter orbifold }

\subsection*{\nmb.{6.1}. Coxeter orbifolds }
A smooth  orbifold $X$ is called a {\it Coxeter orbifold} if for each
local uniformizing system
$(\tilde U_i\subset \mathbb R^n, G_i\subset O(n), \ph_i)$ the
group $G_i$ is a finite linear Coxeter group.

\subsection*{\nmb.{6.2}. Example }
Let $M$ be a Coxeter Riemannian manifolds with reflection group $G$.
Then any Weyl chamber is a Coxeter orbifold.
This follows from proposition \nmb!{5.1} and \nmb!{2.11}.

\subsection*{\nmb.{6.3}. Coxeter orbifold as a manifold with corners and its
universal Coxeter equipment}
Let $X$ be a Coxeter orbifold.
Let $(\tilde U_i, G_i,\ph_i)_{i\in I}$ be an atlas of local uniformizing
systems on $X$ such that $(U_i)$ is an open cover of $X$.
Then $\tilde U_i\subset \mathbb R^n$ is an open neighborhood of 0 which is
invariant under the Coxeter group $G_i$. Thus the orbit space $\tilde
U_i/G_i$ is an open neighborhood of 0 in a linear Weyl chamber of the
group $G_i$. The (equivariant) chart changings $\ph_{ij}$ induce smooth
chart changings between open subsets of $\tilde U_j/G_j$
and $\tilde U_i/G_i$. These respect the
indices of corners (see \nmb!{4.1}).
Thus they describe a smooth atlas for the structure of a
manifold with corners on $X$. So walls and faces are defined and to each
wall $W$ one can associate a generator $s(W)$ of the Coxeter system
with the following property: If $W\cap U_i\ne\emptyset$ for a local
uniformizing system $(\tilde U_i, G_i,\ph_i)$, then $s(W)$ equals the
generator of $G_i$ which is given by the reflection in the wall
$\ph_i\i(W)\subset \tilde U_i$.
Then $(s(W)s(W'))^{n(W,W')}=1$ if
$\ph_i\i(W),\ph_i\i(W')\ne\emptyset$ the generators corresponding to them
in $G_i$ satisfy the same relation.

\begin{proclaim}{\nmb.{6.4}.  Theorem}
Any Coxeter orbifold is the Weyl chamber of a Riemannian Coxeter manifold.
\end{proclaim}

\begin{demo}{Proof}
This follows from \nmb!{6.3} and \nmb!{4.7}.
\qed\end{demo}

\begin{proclaim}{\nmb.{6.5}.  Corollary}
Any Coxeter orbifold is good in the sense of Thurston \cit!{31}. 
\end{proclaim}

\subsection*{\nmb.{6.6}.  Coxeter orbifold structures on a simplex}
Let $\Delta_n$ be the standard $n$-simplex with vertices $0,1, \dots ,
n$. If $s$ is a Coxeter $(G,S)$-equipment of $\Delta_n$, then
there exist a Riemannian metric $\gamma$ on $\Delta_n$, such that
$(\Delta_n, \gamma)$ is a Riemannian chamber and the equipment $s$
is admissible. We denote by $M =\mathcal U(G,\Delta_n, s, \gamma)$ the
associated Coxeter $G$-manifold. It is simply connected. The
homeomorphism $M/G \cong \Delta_n$ define on $\Delta_n$ a
structure of Coxeter orbifold, with the universal covering
manifold $M$, which depends only on the equipment $s$, by \nmb!{4.8}.
Hence, a description of Coxeter orbifold structures on $\Delta_n$
and also Coxeter $G$-manifolds  with the orbit space $\Delta_n$
up to a $G$-diffeomorphism  reduces to a description of Coxeter
equipments of $\Delta$.

For any finite Coxeter group $G$ with the generators $S =\{s_0, \dots, s_n \}$
there exist a unique natural equipment such that the
wall $W_i = (0,1,\dots , \widehat i, \dots , n)$ corresponds to
$s_i$ for $i = 0, \dots, n$. The corresponding Coxeter manifold $M$
is  the sphere $S^n$ with the natural action of $G$ induced by
the standard representation of $G$ in $\mathbb R^{n+1}$.
Let now $G$ be an infinite Coxeter group with system of
generators $S$. There exists a Coxeter $(G,S)$-equipment of
$\Delta_n$ if and only if $|S|=n+1$ and the Coxeter subgroup
generated by $S\setminus \{s\}$ is finite for any $s \in S$.  In
term of the Coxeter diagram $\Gamma$ of the group $G$, this means that
all connected components of $\Gamma$ with exception of one
component $\Gamma_0$ correspond to finite Coxeter groups, and the
component $\Gamma_0$ corresponds to  an infinite Coxeter group,
but after deleting any node it become a Coxeter diagram of a
finite Coxeter  group. One can easily check that such a Coxeter
diagram $\Gamma_0$ is  either a connected parabolic Coxeter
diagram (i.e., extended Dynkin diagram of a simple Lie algebra) or
one of the following diagrams:

\centerline{
\begin{picture}(35,36)\put(3,18){\makebox(0,0){$\bullet$}}
                         \put(3,18){\line(2,1){20}}
                         \put(13,26){\makebox(0,0){$r$}}
                         \put(3,18){\line(2,-1){20}}
                         \put(13,8){\makebox(0,0){$p$}}
                         \put(23,28){\makebox(0,0){$\bullet$}}
                         \put(23,28){\line(0,-1){20}}
                         \put(27,18){\makebox(0,0){$q$}}
                         \put(23,8){\makebox(0,0){$\bullet$}}
\end{picture}
\qquad
\begin{picture}(76,12)(0,-15)\put(3,2.7){\makebox(0,0){$\bullet$}}
                         \put(3,4.5){\line(1,0){20}}
                         \put(3,3){\line(1,0){20}}
                         \put(3,1.5){\line(1,0){20}}
                      \put(23,2.7){\makebox(0,0){$\bullet$}}
                         \put(23,3){\line(1,0){20}}
                      \put(43,2.7){\makebox(0,0){$\bullet$}}
                         \put(43,3){\line(1,0){20}}
                      \put(63,2.7){\makebox(0,0){$\bullet$}}
\end{picture}
\qquad
\begin{picture}(76,12)(0,-15)\put(3,2.7){\makebox(0,0){$\bullet$}}
                         \put(3,3){\line(1,0){20}}
                      \put(23,2.7){\makebox(0,0){$\bullet$}}
                         \put(23,4.5){\line(1,0){20}}
                         \put(23,1.5){\line(1,0){20}}
                         \put(23,3){\line(1,0){20}}
                      \put(43,2.7){\makebox(0,0){$\bullet$}}
                         \put(43,3){\line(1,0){20}}
                      \put(63,2.7){\makebox(0,0){$\bullet$}}
\end{picture}
}

\bigskip
\centerline{
\begin{picture}(36,36)(0,-5)\put(3,2.7){\makebox(0,0){$\bullet$}}
                         \put(3,4.5){\line(1,0){20}}
                         \put(3,3){\line(1,0){20}}
                         \put(3,1.5){\line(1,0){20}}
                      \put(23,2.7){\makebox(0,0){$\bullet$}}
                         \put(23,3){\line(0,1){20}}
                      \put(23,22.7){\makebox(0,0){$\bullet$}}
                         \put(23,22.7){\line(-1,0){20}}
                      \put(3,22.7){\makebox(0,0){$\bullet$}}
                       \put(3,22.7){\line(0,-1){20}}
\end{picture}
\qquad
\begin{picture}(76,24)(0,-15)\put(3,2.7){\makebox(0,0){$\bullet$}}
                         \put(3,4.5){\line(1,0){20}}
                         \put(3,3){\line(1,0){20}}
                         \put(3,1.5){\line(1,0){20}}
                      \put(23,2.7){\makebox(0,0){$\bullet$}}
                         \put(23,3){\line(1,0){20}}
                      \put(43,2.7){\makebox(0,0){$\bullet$}}
                         \put(43,3){\line(2,1){20}}
                         \put(43,3){\line(2,-1){20}}
                      \put(63,12.7){\makebox(0,0){$\bullet$}}
                      \put(63,-7.3){\makebox(0,0){$\bullet$}}
\end{picture}
\qquad
\begin{picture}(96,12)(0,-15)\put(3,2.7){\makebox(0,0){$\bullet$}}
                         \put(3,4.5){\line(1,0){20}}
                         \put(3,3){\line(1,0){20}}
                         \put(3,1.5){\line(1,0){20}}
                      \put(23,2.7){\makebox(0,0){$\bullet$}}
                         \put(23,3){\line(1,0){20}}
                      \put(43,2.7){\makebox(0,0){$\bullet$}}
                         \put(43,3){\line(1,0){20}}
                         \put(63,2.7){\makebox(0,0){$\bullet$}}
                         \put(63,3){\line(1,0){20}}
                         \put(83,2.7){\makebox(0,0){$\bullet$}}
\end{picture}
}

An  interesting question is to classify such equipments for other
polyhedra, e.g.,  a cube, prism etc. 

May be, it is possible to
construct a  non trivial example of compact Coxeter  manifolds,
e.g. simply connected 3-manifolds.

\nocite{*}\bibliographystyle{amsplain}

\begin{thebibliography}{10}

\bibitem{\cit0{1}}
Dmitri Alekseevky, Andreas Kriegl, Mark Losik, Michor, and Peter W.,
  \emph{Lifting smooth curves over invariants for representations of compact
  {L}ie groups}, Transformation groups \textbf{5} (2000), no.~2, 103--110.
  arXiv:math.DG/9801029.

\bibitem{\cit0{2}}
\bysame, \emph{The {R}iemannian geometry of orbit spaces. {T}he metric,
  geodesics, and integrable systems}, Publ. Math. \textbf{62} (2003), no.~3-4,
  arXiv:math.DG/0102159.

\bibitem{\cit0{3}}
E.~Bierstone, \emph{Lifting isotopies from orbit spaces}, Topology \textbf{14}
  (1975), 245--252.

\bibitem{\cit0{4}}
N.~Bourbaki, \emph{Groupes et alg\`ebres de {L}ie}, Hermann, Paris, 1968, Ch.
  4--6.

\bibitem{\cit0{5}}
L.~Br\"ocker, \emph{Semialgebraische {G}eometrie}, Jber. d. Dt. Math.-Verein.
  \textbf{97} (1995), 130--156.

\bibitem{\cit0{6}}
Yu. Burago, M.~Gromov, and G.~Perel'man, \emph{A. {D}. {A}lexandrov spaces with
  curvature bounded below}, Russ. Math. Surv. \textbf{47} (1992), 1--58.

\bibitem{\cit0{7}}
Ruth Charney and Michael~W. Davis, \emph{Singular metrics of nonpositive
  curvature on branched covers of {R}iemannian manifolds}, Am. J. Math.
  \textbf{115} (1993), 929--1009.

\bibitem{\cit0{8}}
\bysame, \emph{The {E}uler characteristic of a nonpositively curved, piecewise
  {E}uclidean manifold}, Pac. J. Math. \textbf{171} (1995), 117--137.

\bibitem{\cit0{9}}
\bysame, \emph{The $k(\pi, 1)$-problem for hyperplane complements associated to
  infinite reflection groups}, J. Am. Math. Soc. \textbf{8} (1995), 597--627.

\bibitem{\cit0{10}}
Yuan-Jen Chiang, \emph{Spectral geometry of $v$-manifolds and its application
  to harmonic maps}, Differential geometry. {P}art 1: {P}artial differential
  equations on manifolds, Proc. Symp. Pure Math. 54, Part 1, 1993, pp.~93--99.

\bibitem{\cit0{11}}
H.S.M. Coxeter, \emph{Discrete groups generated by reflections}, Ann. Math.
  \textbf{35} (1934), 388--621.

\bibitem{\cit0{12}}
Michael~W. Davis, \emph{Groups generated by reflections and aspherical
  manifolds not covered by {E}uclidean space}, Ann. Math. \textbf{117} (1983),
  293--324.

\bibitem{\cit0{13}}
\bysame, \emph{Coxeter groups and aspherical manifolds}, Algebraic {T}opology,
  {A}arhus 1982, Lect. Notes Math. 1051, Springer-Verlag, 1984, pp.~197--221.

\bibitem{\cit0{14}}
\bysame, \emph{A hyperbolic 4-manifold}, Proc. Am. Math. Soc. \textbf{93}
  (1985), 325--328.

\bibitem{\cit0{15}}
\bysame, \emph{The homology of a space on which a reflection group acts}, Duke
  Math. J. \textbf{55} (1987), 97--104.

\bibitem{\cit0{16}}
\bysame, \emph{Some aspherical manifolds}, Duke Math. J. \textbf{55} (1987),
  105--139.

\bibitem{\cit0{17}}
\bysame, \emph{Convex polytopes, {C}oxeter manifolds, and torus actions}, Duke
  Math. J. \textbf{62} (1991), 417--455.

\bibitem{\cit0{18}}
A.~Douady and L.~H\'erault, \emph{Arrondissement des {V}ari\'et\'es \`a coins},
  Comm. Math. Helv. \textbf{48} (1973), 484--491, Appendice \`a: A. Borel, J.P.
  Serre: Corners and arithmetic groups.

\bibitem{\cit0{19}}
J.~E. Humphreys, \emph{Reflection groups and {C}oxeter groups}, Cambridge
  University Press, Cambridge, 1990, 1992, Cambridge studies in advanced
  mathematics 29.

\bibitem{\cit0{20}}
N.~Iwahori, \emph{On discrete reflection groups on symmetric {R}iemannian
  manifolds}, pp.~57--62, Nippon Hyoronsha, Tokyo, 1966 w, MR~36\#830,
  Z~146.43501.

\bibitem{\cit0{21}}
R.~Kane, \emph{Reflection groups and invariant theory}, Springer-Verlag, 2001,
  CMS books in Math. 5.

\bibitem{\cit0{22}}
S~{\L}ojasiewicz, \emph{Triangulation of semi-analytic sets}, Ann. Sc. Norm.
  Sup. Pisa, ser. III \textbf{18} (1964), 449--473.

\bibitem{\cit0{23}}
M.V. Losik, \emph{Lifts of diffeomorphisms of orbit spaces for representations
  of compact {L}ie groups}, Geom. Dedicata \textbf{88} (2001), 21--36.

\bibitem{\cit0{24}}
M.A. Michailowa, \emph{Factor space under the action of a finite group
  generated by pseudo reflections}, Izwest. Acad. Nauk USSR, Ser. Math.
  \textbf{48} (1984), 104--126.

\bibitem{\cit0{25}}
Peter~W. Michor, \emph{Manifolds of differentiable mappings}, Shiva Publ.,
  Orpington, 1980, Shiva Mathematics Series 3 MR~83g:58009, ZB~433.58001.

\bibitem{\cit0{26}}
P.W. Michor, \emph{Isometric actions of {L}ie groups and invariants}, Lecture
  course at the University of Vienna, 1996/97, {\tt
  www.univie.ac.at/\~{}michor/tgbook.ps}

\bibitem{\cit0{27}}
Satake, \emph{Gauss-bonnet theorem for $v$-manifolds}, J. Math. Soc. Japan
  \textbf{9} (1957), 464--492.

\bibitem{\cit0{28}}
G.~W. Schwarz, \emph{Smooth functions invariant under the action of a compact
  {L}ie group}, Topology \textbf{14} (1975), 63--68.

\bibitem{\cit0{29}}
\bysame, \emph{Lifting smooth homotopies of orbit spaces}, Publ. Math. IHES
  \textbf{51} (37--136), 1980.

\bibitem{\cit0{30}}
J.~Szenthe, \emph{A significant property of real hyperbolic spaces},
  pp.~153--161, Babe{\c s}-Bolyai University, Cluj-Napoca, 1993, Preprint, 93-2
  96b:53063.

\bibitem{\cit0{31}}
William~P. Thurston, \emph{Three-dimensional geometry and topology. {V}ol. 1.,
  silvio levy, ed.}, Princeton University Press, Princeton, 1997, Princeton
  Mathematical Series 35.

\bibitem{\cit0{32}}
E.~B. Vinberg, \emph{Discrete linear groups that are generated by reflections},
  Izv. Akad. Nauk SSSR Ser. Mat. \textbf{35} (1971), 1072--1112, English
  translation, Math. USSR, Izv., {\bf 5} (1971), 1083-1119.

\bibitem{\cit0{33}}
E.B. Vinberg, \emph{Discrete reflection groups in {L}obachevsky spaces}, Proc.
  {I}nt. {C}ongr. {M}ath., {W}arszawa 1983, {V}ol. 1, 1984, pp.~593--601.

\bibitem{\cit0{34}}
\bysame, \emph{Hyperbolic reflection groups}, Russ. Math. Surv. \textbf{40}
  (1985), 31--75.

\bibitem{\cit0{35}}
E.B. Vinberg and O.V. Shvartsman, \emph{Discrete groups of motions of spaces of
  constant curvature}, Geometry. {I}I: {S}paces of constant curvature,
  {E}ncycl. {M}ath. {S}ci. 29 (Berlin), Springer-Verlag, 1993, pp.~139--248.

\bibitem{\cit0{36}}
H.~Whitney, \emph{Elementary structure of real algebraic varieties}, Ann. Math.
  \textbf{66} (1957), 545--556.

\end{thebibliography}
\providecommand{\bysame}{\leavevmode\hbox to3em{\hrulefill}\thinspace}
\providecommand{\MR}{\relax\ifhmode\unskip\space\fi MR }
\providecommand{\MRhref}[2]{%
  \href{http://www.ams.org/mathscinet-getitem?mr=#1}{#2}
}
\providecommand{\href}[2]{#2}

\end{document}